\documentclass[12pt]{article}
\usepackage{amssymb,latexsym,amsmath}

\newtheorem{corollary}{Corollary}
\newtheorem{lemma}[corollary]{Lemma}
\newtheorem{definition}[corollary]{Definition}
\newtheorem{proposition}[corollary]{Proposition}
\newtheorem{theorem}[corollary]{Theorem}

\begin{document}


\title{The Fredholm index of a pair of commuting operators}
\author{Xiang Fang
\footnote{Partially supported by National Science Foundation Grant
DMS 0400509}
}
\date{}

\maketitle

\def\cc{\mathbb{C}}
\def\zz{\mathbb{Z}}
\def\nn{\mathbb{N}}
\def\rr{\mathbb{R}}
\def\qq{\mathbb{Q}}
\def\dd{\mathbb{D}}
\def\smult{s_{\_}mul(T)}
\def\smul#1{s_{\_}mul(\emph{\emph{\emph{}}}{#1})}
\def\bsmult{b.s._{\_}mul(T)}
\def\bsmul#1{b.s._{\_}mul({#1})}
\def\submo{\mathcal{M}}
\def\submoperp{{\submo}^{\perp}}

\def\tinbh{T \in B(M)}

\def\bigno{\bigskip \noindent}

\begin{abstract}
This paper concerns   Fredholm theory in several variables, and
  its applications to   Hilbert spaces of analytic
functions.  One feature is       the introduction
  of ideas from   commutative algebra  to operator theory.

\medskip

Specifically, we introduce a method to calculate the   Fredholm
index of a pair of commuting operators.  To achieve this, we
define and study the Hilbert space analogs of Samuel
multiplicities in commutative algebra.

\medskip

 Then the theory is applied to the   symmetric Fock space.
%
In particular,  our results imply a   satisfactory answer to
Arveson's program on developing a  Fredholm theory for pure
$d$-contractions  when $d=2$, including both the Fredholmness
problem and the calculation of indices.
 We also show that  Arveson's curvature invariant is in fact always equal to the Samuel
multiplicity for an arbitrary pure d-contraction with finite
defect rank.
 It follows that the curvature is a similarity
invariant.
\end{abstract}

\newpage
\tableofcontents







\newpage
\setcounter{section}{-1}
%
\section{Introduction}

\subsection{Background and motivation}
The idea of developing    Fredholm theory, beyond that of a single
operator, has been around for many years. In particular, a theory
  based on J. Taylor's Koszul complex approach to multivariable
spectral theory \cite{Taylor1}, \cite{Taylor2} seems particularly
appealing,  and can be formulated  as the index theory of  an
abstract Dirac operator \cite{Dirac}, which can be regarded as an
abstraction of the local model of Dirac operators (\cite{Higson},
page 310) in Riemannian geometry. But currently there are
essentially no effective tools available to calculate the Fredholm
index in several variables. In this paper we will introduce a
general method to calculate the index in the two variable case,
which can be regarded as a Hilbert space version of a classical
theorem of J.-P. Serre in local algebra (\cite{Serre}, page 57,
Theorem 1).
A novelty of our approach is the extensive use of ideas from
commutative algebra in a Hilbert space setting.

Among the applications to operator theory we mention two results
next. Theorem \ref{T:subFock} implies a satisfactory answer to
Arveson's program \cite{Arveson-problem} on developing a Fredholm
theory for pure $d$-contractions when $d=2$, including both the
Fredholmness problem and the calculation of indices. It is
noteworthy   that the Fredholmness problem is usually quite
subtle, and previous results are scarce \cite{Arveson-psum},
\cite{Arveson-problem}.
Theorem  \ref{T:cur} shows that  Arveson's curvature invariant is
in fact always equal to the Samuel multiplicity, a notion borrowed
from algebra, for an arbitrary pure d-contraction with finite
defect rank. This provides an answer to Arveson's question
\cite{Dirac}, \cite{Arveson-psum} on  how to express the curvature
invariant in terms of other invariants which are determined by the
d-contraction directly, and can immediately imply that the
curvature is an integer. From Theorem \ref{T:cur} it follows that
the curvature is a similarity invariant, which is previously
unknown.

\medskip

Next we discuss some background information and motivations. In
the past
  much of
the   effort in higher dimensional Fredholm theory  was devoted to
the study of general Fredholm complexes, especially their
stability under various
   perturbations.
%
See Ambrozie-Vasilescu \cite{Ambrozie}, Curto \cite{curtothesis},
\cite{Curto}, Eschmeier-Putinar \cite{E-Putinar}, Segal
\cite{Segal}, Vasilescu \cite{Vasilescu}, and the references
therein (most of these   contain   extensive bibliographies).
What is more relevant to this paper is the study of
 numerical invariants. Along this line, Carey-Pincus \cite{Carey-Pincus}, Levy
\cite{Levy1}, \cite{Levy2}, Putinar \cite{Putinar} and the last
chapter of the book of Eschmeier-Putinar \cite{E-Putinar},
establish many connections between multivariable Fredholm theory
and various areas in mathematics, such as K-theory, index theory,
geometric measure theory, sheaf theory, and even cyclic cohomology
in noncommutative geometry. Also, see Arveson's expository paper
\cite{Arveson-problem} for a recent development along a different
line.

However, as far as the calculation of multivariable Fredholm
indices is concerned, still no effective machineries are
available. Although there are many evidences showing that such
calculations are valuable and should have close connections with
index theory of Dirac operators in geometry \cite{Levy},
\cite{Segal} and function theory over various domains in $\cc^n$
\cite{E-Putinar}.

\medskip

 What we want to emphasize in this paper is
that Koszul complexes are not just general Fredholm complexes:
there exist many commutative algebraic features, such as tensorial
properties, which are particular to Koszul complexes, and not
enjoyed by general Fredholm complexes. These features play
significant roles in many areas of commutative algebra, but  have
not been exploited in operator theory.
%
%
%
One of the contributing factors to this situation might be that
Koszul complexes over Hilbert spaces are somehow difficult to
handle from both algebraic and analytic aspects: Hilbertian Koszul
complexes lack basic properties such as
  Noetherian conditions to be algebraically  amenable, while their
  algebraic structure is not familiar to many
  analysts.

Our motivation is the following observation.  J. Taylor's work
shows that in order to define fundamental concepts such as
invertibility in multivariable operator theory, objects
(\emph{i.e.}, Koszul complexes) with a rich algebraic theory
necessarily enter   the picture; it suggests that  a well
developed multivariable operator theory will be inextricably
interwoven with homological and commutative algebraic methods. The
Hilbert module program of Douglas-Paulsen \cite{DP} provides such
a framework.

On the other hand, J. Taylor's framework is well suited for
Fredholm theory \cite{curtothesis}.
This leads us to  think about :

\medskip

(A) how to formulate Fredholm theory from an algebraic viewpoint;
and

(B) how to   connect existing parts of operator theory   to this
algebraic viewpoint.

\medskip

What appears to be plausible to us is that, for (A)  Fredholm
theory can be studied as a Koszul homology, defined on topological
modules (or, Hilbert modules \cite{DP}), and we seek to develop
methods to calculate this (topological) homology theory;
for (B), we can interpret natural operator theoretic invariants,
such as the kernel $ker(T)$ and cokernel $H/TH$ of an operator $T
\in B(H)$, in terms of homological invariants:
 they are just the
two homology groups of the length one Koszul complex $0 \to H
\xrightarrow{T} H \to 0$.

More generally, a systematic way to introduce homological and
commutative algebraic methods into operator theory is described in
 \cite{proposal}.
   For a single
operator, this program turns out to be fruitful  and leads to
  results in operator theory \cite{Fred} ,
\cite{Dirichlet}.
 This paper represents some progress in the more than
one variable case.

\bigskip

 Thus the purpose of this paper is twofold.
 Firstly, we show that many
algebraic methods carry over to the study of Koszul complexes over
Hilbert modules; in the meanwhile, these algebraic considerations
lead to new phenomena, which are unseen for algebraic modules (see
Subsection \ref{subS:index}).
Secondly, and perhaps more importantly, we show that these purely
algebraic machineries interact well with function theory and
operator theory, and indeed lead to operator theoretic results.

 \bigskip

Since the paper treats both algebraic theory and function theory,
sometimes we have provided more details than necessary, for the
benefit  of readers with different backgrounds. In particular,
operator theorists might be interested in more references in
algebra: Eisenbud's book \cite{Eisenbud} is a standard resource
for Koszul complexes, Hilbert polynomials and Samuel
multiplicities; the expository paper \cite{Auslander} of
Auslander-Buchbaum provides a nice treatment of Serre's
multiplicity theory; various ideas in Serre's book \cite{Serre},
P. Roberts' book \cite{Roberts}, and the papers of C. Lech
\cite{Lech1}, \cite{Lech2} have been very helpful. Also
Balcerzyk-Jozefiak's book \cite{Balcerzyk} contains another proof
of Serre's theorem on the Euler characteristic and Samuel
multiplicity.

\subsection{Organization}
The body of the paper contains four sections.

  Section 1   develops the necessary algebraic theory for Koszul
complexes over Hilbert spaces, which is needed in subsequent
sections.
 In particular, it   introduces
 a set of refined invariants to calculate the Fredholm index of
any Fredholm pair of commuting operators (\ref{subS:index}).
Moreover, it contains  induction arguments (\ref{subS:induction}),
  estimates of Hilbert functions (\ref{subS:estimate}), and
  examples.

Section 2 contains the applications to the Fredholm theory on the
$L$-valued symmetric Fock space in two variables (Theorem
\ref{T:subFock}). Here $L$ is a
  separable Hilbert space,  not necessarily finite dimensional. In particular, we
completely characterize the Fredholmness of the pair of
multiplication by coordinate functions   on invariant subspaces.
Then the Fredholm index is calculated in terms of the Samuel
multiplicity.

Section 3 establishes a version of the additivity of Hilbert
polynomials (Theorem \ref{T:add}).
It follows a monotonicity property of Samuel multiplicities
(Corollary \ref{C:mono}), which can be regarded as a higher
dimensional generalization of the codimension-one property for
invariant subspaces of the Hardy or Dirichlet space over the unit
disc.

The study of this additivity property is intimately related to the
study of the curvature invariant introduced by Arveson
\cite{curvature}. We show that the curvature invariant is always
equal to the Samuel multiplicity for an arbitrary pure
d-contraction with finite defect rank (Theorem \ref{T:cur}). As a
consequence, the curvature invariant is shown to be preserved
under similarity relations (Corollary \ref{C:sim}).

Section 4 provides one more proof of the index formula established
in Theorem \ref{T:main}, Part (ii), and some comments concerning
further studies.

\subsection{The one variable case}
The one variable case can serve as a good illustration of the
basic ideas, which we shall explore. Moreover, this paper makes
essential use of the one variable theory developed in \cite{Fred}.
Because we shall apply results in \cite{Fred} (mainly Theorem 2
and its proof) in a non-Hilbert-space setting, which is more
general than what is actually proved in that paper,  not only
shall we  use the theorems there, but also the arguments to prove
them.
 Since  it does not seem worthwhile to give a separate treatment of
this non-Hilbert-space case, we advise   readers unfamiliar with
\cite{Fred} to go through the proof of Theorem 2 there.

Recall that, for a single Fredholm operator $T$, its index is
defined to be
\begin{equation*}
indexs(T)=dim \ ker(T) - dim \ H/TH.
\end{equation*}
 Then $index(T)$ has
remarkable stability, while, individually, $ker(T)$ and $H/TH$ are
sensitive to perturbations. The starting point is the simple fact
\begin{equation}\label{E:dynamics}
index(T^k)=k \cdot index(T),
\end{equation}
or, equivalently,
\begin{equation}\label{E:dynamics2}
index(T) =\frac{dim \ ker(T^k)}{k} - \frac{ dim \ H/T^kH} {k}.
\end{equation}
 So, instead of $ker(T)$ and $H/TH$, we
look at the asymptotic behavior of $ker(T^k)$ and $H/T^kH$ as $k
\to \infty$. At this point, experts can immediately recognize that
this will connect  operator theory with Hilbert polynomials in
algebraic geometry, and the asymptotic behavior of Grothendieck's
local cohomology modules \cite{Brodmann}, which is a currently
active research area in commutative algebra. These connections
seem promising and are not fully understood.

When $k>>0$, one has
\begin{eqnarray*}
dim \ ker(T^k) & = & ak+b \\
  dim \ H/ T^kH &=& ck+d
\end{eqnarray*}
for some integers $a$, $b$, $c$, and $d$. These equalities  might
not be obvious, but they are simple consequences of the existence
of Hilbert polynomials, or they can be directly proved by showing
that, the dimension of the difference $T^kH \ominus T^{k+1}H$
stablizes, as $k \to \infty$. By Equation (\ref{E:dynamics2})
\begin{displaymath}
index(T)=a-c.
\end{displaymath}

Formally, $a$ and $c$ are analogous to  Samuel multiplicities in
algebra. Hence one can ask a variety of questions from an
algebraic viewpoint.  But in the case of operator theory, we feel
that the real significance of the integers $a$ and $c$   lies in
the fact that, they are defined on a module with a topology, that
is, the Hilbert space $H$ as a $\cc[z]$-module, hence they might
have some geometric meaning; when compared with $dim \ ker(T)$ and
$dim \ H/TH$, $a$ and $c$ have more stability, and the stability
might lead to some geometric interpretation of $a$ and $c$. This
is indeed true, and has lead to  results in operator theory
\cite{Fred}, \cite{Dirichlet}.

\bigskip

\section{Algebraic theory: Hilbertian Koszul complexes}
\subsection{Fredholm index and Samuel multiplicity}\label{subS:index}
Given a tuple of commuting operators $T=(T_1, \cdots, T_d)$,
acting on a Hilbert space $H$, we can endow $H$ with an
$A=\cc[z_1, \cdots, z_d]$-module structure by \cite{DP}
\begin{displaymath}
(p(z_1, \cdots, z_d), \xi) \in A \times H \to p(T_1, \cdots,
T_d)\xi \in H.
\end{displaymath}
For any $x\in A$, we denote by
\begin{displaymath}
K(x)= 0 \to A \xrightarrow{x} A \to 0
\end{displaymath}
 the Koszul complex
associated with $x$. For $x_1, \cdots, x_r \in A$, we define the
Koszul complex
\begin{displaymath}
K(x_1, \cdots, x_r)=K(x_1) \otimes_A K(x_2) \otimes_A \cdots
\otimes_A K(x_r).
\end{displaymath}
See \cite{Eisenbud} for basic definitions about tensor products of
complexes. For any $A$-module $M$, we define $K(x_1, \cdots, x_r;
M)=M\otimes_A K(x_1, \cdots, x_r)$.
 In the case of a Hilbert space $H$,
whose $A$-module structure is induced by a commuting operator
tuple $T=(T_1, \cdots, T_d)$, we also write $K(T_1, \cdots, T_d;
H) = K(z_1, \cdots, z_d; H)$. Let $H_iK(\  \cdots \ ; M)$ denote
the $i$th homology group of a Koszul complex $K(\  \cdots \ ; M)$.

Then   comes J. Taylor's definition of invertibility: for a single
operator $T$ we say that $T$ is invertible if and only if
\begin{equation*}
H\otimes_{\cc[z]}K(z)
\end{equation*}
 is exact, equivalently, $0\to H\xrightarrow{T} H \to 0$ is exact; for a commuting tuple $T=(T_1,
\cdots, T_d)$,   $T$ is invertible if and only if
\begin{equation*}
H\otimes_A K(z_1)\otimes_A \cdots \otimes_A K(z_d)
\end{equation*}
 is exact.
 We say that $(T_1, \cdots, T_d)$ is Fredholm if each
$H_iK(T_1, \cdots, T_d; H)$ is finite dimensional as a vector
space; then the Fredholm index of the tuple $T$ is defined as
\begin{displaymath}
index(T)=(-1)^d \  \chi(K(T; H))=\sum_{i=0}^d (-1)^{d-i} \ dim \
H_iK(T_1, \cdots, T_d; H).
\end{displaymath}

  Our definition of the Koszul complex and
Taylor invertibility has a strong algebraic flavor, and looks
different from some of the past references in multivariable
operator theory, say \cite{Dirac}, \cite{Curto}, \cite{E-Putinar},
\cite{Vasilescu}. They are of course equivalent. Our definition
will allow one to work with algebraic properties more naturally.

 In the case of a pair, note that $H_2K(T; H)=ker(T_1)\cap
ker(T_2)$, and $H_0K(T;H)=H/(T_1H+T_2H)$. In particular,
$T_1H+T_2H$ is automatically closed when $T$ is Fredholm.

\bigskip
Because of the following Proposition \ref{P:xx'}, we shall be able
to adopt a strategy similar to that used for the study of a single
operator.

\begin{proposition}\label{P:xx'}
Let $(T_1, T'_1, T_2, \cdots, T_d)$ be a $(d+1)$ tuple of
commuting operators, acting on a Hilbert space $H$.  If any two of
the following three d-tuples $(T_1, T_2, \cdots, T_d)$, $(T'_1,
T_2, \cdots, T_d)$, and  $(T_1 T'_1, T_2, \cdots, T_d)$ are
Fredholm, then so is the third one.  When all are Fredholm, we
have
\begin{displaymath}
index(T_1  T'_1, T_2, \cdots, T_d) = index(T_1, T_2, \cdots, T_d)
+ index(T'_1, T_2, \cdots, T_d).
\end{displaymath}
\end{proposition}

\begin{corollary}\label{C:dyna}
For any  d-tuple $T=(T_1, \cdots, T_d)$ of commuting operators, if
$T$ is Fredholm, then
\begin{equation*}
index(T_1^{n_1}, \cdots, T_d^{n_d})=n_1\cdots n_d \ index(T_1,
\cdots, T_d)
\end{equation*}
for all $n_1, \cdots, n_d  \in \nn$.
\end{corollary}

\noindent \textbf{Proof of Proposition \ref{P:xx'}.} The proof is
purely algebraic, and we adopt the arguments in \cite{Roberts}
page 97. We also call the readers attention to Putinar's
\cite{Putinar-alge} where many of the basic algebraic properties
of a Koszul complex were established.

\bigskip

 Consider $H$ as a Hilbert module over   the
polynomial ring $A$ in $(d+1)$ variables, with module actions
given by $(T_1, T'_1, T_2, \cdots, T_d)$.

 First we remark that, if we also consider $H$ as a Hilbert module
over the polynomial ring in $d$ variables, with module actions
given by $(T_1, T_2, \cdots, T_d)$, then there are two natural
ways to associate a Koszul complex to the tuple $(T_1, T_2,
\cdots, T_d)$. The resulting Koszul complexes $K(T_1, \cdots, T_d;
H)$ are in fact independent of the base ring and hence produce the
same Fredholm index.

Let $E_.=0\to E_1 \xrightarrow{d_1} E_0 \to 0$ be the complex
given by $E_1=E_0=A\oplus A$, and $d_1=\left(
\begin{array}{cc}z_1 & 1 \\ 0 & z_2 \end{array}\right)$.
Then we have an exact sequence
\begin{equation}\label{E:x1x2}
0\to K(z_1; A) \to E_. \to K(z_2; A) \to 0
\end{equation}
by embedding $K(z_1; A)$ into the first copy of $A$ of $E_1$ and
$E_0$.

Next, we relate $E_.$ to $K(z_1z_2; A)$ by
%
%
the following exact sequence
\begin{equation}\label{E:x1x2again}
0\to K(z_1z_2; A) \xrightarrow{f_.} E_. \xrightarrow{g_.}
K(\textbf{1}, A) \to 0.
\end{equation}
Here $f_1=(-1, z_1)$, $f_0=(0, 1)$, $g_1=(z_1, 1)$, and $g_0=(1,
0)$.

 Now apply $H\otimes_A \cdot \otimes_A K(z_3, \cdots,
z_{d+1})$ to sequence (\ref{E:x1x2}) and (\ref{E:x1x2again}).
Recall that tensoring with a complex of free modules preserves
exactness. Moreover,  $\chi(K(1, z_3, \cdots, z_{d+1}; H))=0$,
since the d-tuple $(1, T_3, \cdots, T_{d+1})$ is invertible. Now
we can complete the proof by using the additivity of the Euler
characteristics. $\Box$

\bigskip

In particular, we have, for any $k \in \nn$,
\begin{equation}\label{E:dyna}
index(T_1^k, T_2^k)=k^2 \ index(T_1, T_2).
\end{equation}
So for a pair $T=(T_1, T_2)$, we define the homological Hilbert
functions
\begin{equation}\label{E:hik}
h_i(k)=dim \ H_i K(T_1^k, T_2^k; H)
\end{equation}
for $i=0, 1, 2$, $k \in \nn$. See \cite{Brodmann}, \cite{Roberts}
for similar definitions
 in algebra. Then
\begin{equation}\label{E:dyna2}
index(T_1, T_2)=\frac{h_2(k)}{k^2} - \frac{h_1(k)}{k^2} +
\frac{h_0(k)}{k^2}.
\end{equation}
Now by analogy with Hilbert polynomials, it is natural to expect
that each $h_i(k)$ is a polynomial in $k$, when $k>>0$, and the
leading coefficients of $h_i$ are the two-variable-analogs of the
numbers $a$ and $c$,   as discussed in the introduction.
Unfortunately, this is not the case, as illustrated by the
following example.

\bigno \textbf{Example} $h_i(k)$ may not be a polynomial when $k
>>0$.

\bigskip

For a similar result in algebra, see \cite{Garcia}. The author
thanks P. Roberts for bringing this paper to his attention.
However, there exists a related polynomial in terms of the length
of local cohomology modules in the graded case (\cite{Brodmann},
page 317, Theorem 17.1.9).

\bigskip

Let $H^2(\mathbb{D})$ be the Hardy space over the unit disc, and
$M_z$ be the multiplication by the coordinate function $z$. Let
\begin{displaymath}
H=H^2(\mathbb{D})\oplus H^2(\mathbb{D}), \quad
T_1=\left(\begin{array}{cc} M_z & 0 \\ 0 & M_z \end{array}
\right), \quad  and \ \
T_2=\left(\begin{array}{cc} 0 & M_z^2 \\ 1 & 0 \end{array}
\right).
\end{displaymath}
Then $T_1$, $T_2$ commute. Observe that
\begin{displaymath}
T_2^{2t}=\left(\begin{array}{cc} M_z^{2t} & 0 \\ 0 & M_z^{2t}
\end{array} \right), \quad  and \ \
T_2^{2t+1}=\left(\begin{array}{cc} 0 & M_z^{2t+2} \\ M_z^{2t} & 0
\end{array}\right).
\end{displaymath}
So
\begin{eqnarray*}
h_0(k)  & = & 2k,  \qquad  \qquad   k  \ even;  \\
 &  =  & 2k-1,   \qquad    k  \  odd,
\end{eqnarray*}
which is not a polynomial when $k>>0$.

Because $h_2(k)=0$ for all $k$, it follows from Equation
(\ref{E:dyna}) that, $h_1(k)$ is also not a polynomial, when
$k>>0$.

The function $h_2(k)$ may not be a polynomial for $k>>0$ follows
from considering $(T_1^*, T_2^*)$.

\bigskip

However, we have
\begin{theorem}\label{T:main}
Let $T=(T_1, T_2)$ be a pair of commuting operators, acting on a
Hilbert space $H$. Assume that $T$ is Fredholm. Then
\begin{itemize}
\item[(i)]  $h_i(k)$ $(i =0, 1, 2)$ may not be a polynomial for
$k>>0$, but the limit
\begin{displaymath}
e_i=e_i(T)=\lim_{k\to \infty} \frac{h_i(k)}{k^2},
\end{displaymath}
exists, and is an integer. In particular, each $e_i$ can be
strictly positive.

\item[(ii)] $index(T_1, T_2)=e_2-e_1+e_0$.

\item[(iii)] $e_i(T_1^t, T_2^s)=t \cdot s \cdot e_i(T_1, T_2)$,
for $t, s \ge 0$, $i=0, 1, 2$.
\end{itemize}
\end{theorem}
\textbf{Proof. }  Parts (i) and (ii) of the theorem can be quickly
reduced to the existence of $e_0$ (Subsection \ref{subS:Lech}), as
follows: First we have $H/(T_1H+T_2H) \cong ker(T_1^*)\cap
ker(T_2^*)$ when $T_1H+T_2H$ is closed. It follows that, $e_2(T_1,
T_2)=e_0(T_1^*, T_2^*)$ exists, if $e_0$ always exists. Then, by
Equation (\ref{E:dyna2}), we conclude that $e_1$ exists, and the
formula in Theorem \ref{T:main} Part (ii) is true.

 Part (iii) follows from Corollary \ref{C:Lech} in Subsection
\ref{subS:estimate}. $\quad \square$

\bigno \textbf{Remark} The additional property of $e_i$ as
displayed in Part (iii)   is well known for the Euler
characteristic in algebra (Proposition \ref{P:xx'}, or
\cite{Roberts}, page 98, Corollary 5.2.4), but the fact that, in
operator theory, each $e_i$ satisfies the formula in Part (iii) is
somehow unexpected.

Also, the fact that, each $e_i$ can be strictly positive, stands
in strong contrast with algebraic results because the algebraic
analogs of $e_i$ $(I\ne 0)$ are always zero under mild conditions.

\bigno \textbf{Example}  Each $e_i$ can be strictly positive.

\bigskip

The following example was communicated to the author by R. Curto
\cite{Curto-example}, and satisfies $e_0=e_2=0$, and $e_1=1$. It
is trivial to see that $e_0=e_1=0$, $e_2=1$, or $e_1=e_2=0$,
$e_0=1$ can happen.

Let $S=\{(s_1, s_2) \in \zz^2, s_1 \ge 0, \ or \ s_2 \ge 0\}$, and
$H$ be the Hilbert space with an orthonormal basis $\{e_{s_1,
s_2}, (s_1, s_2)\in S\}$. Define $T_1$, $T_2$ by $T_1e_{s_1,
s_2}=e_{s_1+1, s_2}$, and $T_2e_{s_1, s_2}=e_{s_1, s_2+1}$. Then
$h_0(k)=h_2(k)=0$, and $h_1(k)=k^2$.

\bigskip

Note that $e_i >0$ for $i\ne 0$ is new to operator theory, hence
there is no corresponding results in algebra. See \cite{Roberts}
page 99 Proposition 5.2.6 for more details in the corresponding
algebraic situation. It in fact represents an essential difference
between algebraic modules and Hilbertian modules. In the algebraic
case, $h_i(k)$ is usually bounded by a polynomial in $k$ with
degree $d-i$ (\cite{Roberts}, page 99). In our case, $d=2$. Hence,
$e_1=e_2=0$, and  they are essentially \emph{invisible} in
algebra. It follows that $\chi(K(T; H))=e_0$, which is a classical
theorem of Serre \cite{Serre}. In this sense our theorem is an
extension of Serre's theorem to Hilbert modules.

   An alternative way to look at the asymptotic
behavior of $h_i(k)$  in algebra is to consider $\lim_{k\to
\infty}\frac{h_i(k)}{k^{d-i}}$. However, the problem of whether or
not these limits exist   seems still open  \cite{Kirby}.

\subsection{Lech's formula}\label{subS:Lech}
 Now we turn to the existence of $e_0$,
which follows immediately from   Lech's formula (Theorem
\ref{T:H-Lech}), together with some basic facts on Hilbert
polynomials which we recall now.

Let $H$ be a Hilbert module over $A=\cc[z_1, \cdots, z_d]$, such
that $dim \ H/\overline{IH} < \infty$. Here $I=(z_1, \cdots,
z_d)\subset A$ is the maximal ideal at the origin, and the bar
denotes Hilbertian closure.
If we form the graded module
\begin{displaymath}
gr(H)=\oplus_{k \ge 0} \frac{\overline{I^kH}}{\overline{I^{k+1}H}}
\end{displaymath}
 over $A$,
then $gr(H)$ is generated by its first component
$H/\overline{IH}$, which is finite dimensional. (Note that $H$ is
usually not generated by $H\ominus \overline{IH}$.) Hence, by
basic results on Hilbert polynomials in algebra \cite{DP},
\cite{Eisenbud},  the function $\phi_H(k)= dim \
H/\overline{I^kH}$, which counts the dimension of the first $k$
components of $gr(H)$, is actually a polynomial of degree at most
$d$, when $k>>0$. Moreover, its leading term (or, the degree d
term) has the form $\frac{e}{d!}k^d$, where $e$ is an integer.
The generating property of the first component implies
\begin{equation}\label{E:gee}
dim \ H/\overline{IH} \ge e.
\end{equation}
We call $e=e(I; H)$ the Samuel multiplicity of $H$ with respect to
$I$. Also note that, in the above discussion, we can replace $I$
by a general ideal $J \subset A$.

\begin {theorem}[Lech's formula]\label{T:H-Lech}
Let $H$ be a Hilbert module over $\cc[z_1, \cdots, z_d]$, such
that $dim \ H/\overline{IH} < \infty$. Let $J_k=(z_1^k, \cdots,
z_d^k)$. Then
\begin{equation*}
\lim_{k\to \infty} \frac{dim \ H/\overline{J_kH}}{k^d}=d! \
\lim_{k\to \infty} \frac{dim \ H/\overline{I^kH}}{k^d}.
\end{equation*}
\end{theorem}

Theorem \ref{T:H-Lech} is apparently a Hilbertian  version of the
following algebraic Lech's formula (Lemma
\ref{L:01algebraicLech}), which plays a key role in our proof of
Theorem \ref{T:H-Lech}. Lemma \ref{L:01algebraicLech} is excerpted
from P. Robert's book \cite{Roberts} (page 101 Theorem
 5.2.8.).
\begin{lemma}\label{L:01algebraicLech}
Let $\alpha=(x_1, \cdots, x_d)$ be an ideal in a local ring $R$,
and let $M$ be a finitely generated $R$-module, such that $M/(x_1,
\cdots, x_d)M$ has finite length. Then
\begin{equation*}
\lim_{k\to \infty} \frac{length(M/(x_1^k, \cdots, x_d^k)M)}{k^d} =
d! \ \lim_{k \to \infty} \frac{length(M/\alpha^kM)}{k^d}.
\end{equation*}
\end{lemma}
Note that when the module $M$ is a $\cc$-vector space, the above
$length$ can be replaced by $dim$. Also, it is a standard fact in
algebra that  the above statement in  local rings implies  a
corresponding version for graded rings/modules.

\bigno \textbf{Proof of Theorem \ref{T:H-Lech}.} By the definition
of $gr(H)$, one has $dim \ H/\overline{I^kH} =dim \
gr(H)/I^kgr(H)$. Apply the algebraic  formula of Lech (Lemma
\ref{L:01algebraicLech}) to the graded module $gr(H)$, we have
\begin{eqnarray*}
d! \ \lim_{k\to \infty} \frac{dim \ H/\overline{I^kH}}{k^d} & = &
d! \
\lim_{k\to \infty} \frac{dim \ gr(H)/ I^k gr(H) }{k^d} \\
 & = & \lim_{k \to \infty} \frac{dim \ gr(H)/ J_k gr(H)}{k^d}.
\end{eqnarray*}
Now consider the natural surjective map, induced by inclusions
\begin{equation}\label{E:gradedback}
gr(H)=\frac{H}{\overline{IH}} \oplus
\frac{\overline{IH}}{\overline{I^2H}} \oplus \cdots \ \to \
\frac{H}{\overline{IH+J_kH}} \oplus
\frac{\overline{IH+J_kH}}{\overline{I^2H+J_kH}} \oplus \cdots,
\end{equation}
whose kernel contains the graded submodule $J_k gr(H)$.

Observe that the latter graded module in (\ref{E:gradedback})
stops after a finite number of steps, since $I^t \subset J_k$ when
$t$ is large enough. Moreover, it gives  rise to a grading of
$H/\overline{J_kH}$. It follows that
\begin{equation}\label{E:wronglabel}
dim \ gr(H)/J_k gr(H) \ge dim \ H/\overline{J_k H},
\end{equation}
 which yields one direction of Theorem \ref{T:H-Lech}
\begin{equation*}
d! \ \lim_{k\to \infty} \frac{dim \ H/\overline{I^kH}}{k^d} \ge
\overline{\lim}_{k\to \infty} \frac{dim \ H/\overline{J_kH}}{k^d}.
\end{equation*}

On the other hand, we have $dim \ H/ \overline{J_kH} \ge e(J_k,
H)$, according to inequality  (\ref{E:gee}). We claim that
\begin{equation*}
e(J_k, H)=e(I^k, H)=k^d \cdot e(I, H).
\end{equation*}
 This claim follows from the
following fact
\begin{equation*}
(I^k)^{t+d} \subset J_k^t \subset (I^k)^t,
\end{equation*}
and the definition of $e(\cdot, H)$.   Thus the other direction of
Theorem \ref{T:H-Lech} follows
\begin{equation*}
\underline{\lim}_{k\to \infty} \frac{dim \ H/\overline{J_kH}}{k^d}
\ge e(I, H). \quad \square
\end{equation*}

Now the proof of Part (i) and (ii) of Theorem \ref{T:main} is
completed. It also follows  that, when $d=2$, the Samuel
multiplicity $e(I, H)$ coincides with   $e_0$, defined in Theorem
\ref{T:main}.

\bigno \textbf{Remark} Now we want to make some
 comparisons between the one and two variable cases. So far,
the idea to calculate the Fredholm index of a pair is parallel to
that of the one variable theory in \cite{Fred}, even though for a
pair the situation is more complicated. But, realistically, the
numbers $e_i$ are still quite difficult to calculate, especially
$e_1$. So we need to develop more techniques (see
\ref{subS:induction}, \ref{subS:estimate}), before we can
effectively apply Theorem \ref{T:main} to the symmetric Fock
space. Also, it is not clear how to obtain a matrix decomposition
for the pair, based on $e_i$, as was done in \cite{Fred}, Theorem
4.

\subsection{Induction arguments}\label{subS:induction}
  In algebra an advantage of
considering Koszul complexes is that, they are well suited for
induction arguments. By this it often means to consider modules of
the form $H/z_1H$, on which the $z_1$-module-action is zero. This
consideration turns out to be helpful in Section 2 of this paper.
Since our main interests lie in the Fredholm index and the Samuel
multiplicity,  next we study these two invariants on $H/T_1H$.

Unlike in Lech's formula, we cannot take the closure of $T_1H$,
which means that, we have to enlarge the Hilbert modules to a
larger category to include non-complete spaces.

\bigskip

Let $H$ be a Hilbert module over $A=\cc[z_1, \cdots, z_d]$, with
module actions given by a d-tuple $(T_1, \cdots, T_d)$.
  Then $H/T_1H$
admits a natural module structures over both $A$ and  $A'=\cc[z_2,
\cdots, z_d]$, where the module action  over $A'$ is given by
$(\tilde{T}_2, \cdots, \tilde{T}_{d})$, which is induced by $(T_2,
\cdots, T_{d})$.
We can similarly consider the Koszul complex $K(z_2, \cdots, z_d;
H/T_1H)=K(\tilde{T}_2, \cdots, \tilde{T}_{d}; H/T_1H)$, and define
the Fredholm index, when all involved homology groups are finite
dimensional vector spaces.

\begin{proposition}\label{P:reduction-index}
If $ker(T_1)=\{0\}$, then
\begin{displaymath}
index(T_1, \cdots, T_d; H) =-index(\tilde{T}_2, \cdots,
\tilde{T}_{d}; H/T_1H).
\end{displaymath}
\end{proposition}
\textbf{Proof.} Observe that
\begin{displaymath}
K(T_1, \cdots, T_d; H)=K(T_1; H) \otimes_A K(z_2, \cdots, z_{d}),
\end{displaymath}
and
\begin{displaymath}
K(\tilde{T}_2, \cdots, \tilde{T}_{d}; H/T_1H)=H/T_1H \otimes_A
K(z_2, \cdots, z_{d}).
\end{displaymath}
So the key is to relate $K(T_1; H)$ to  $H/T_1H$. We begin with
the exact sequence of $A$-modules
\begin{displaymath}
0 \to H \xrightarrow{T_1} H \to H/T_1H \to 0.
\end{displaymath}
Then extend it to an exact sequence of complexes
\begin{displaymath}
0 \to K(\textbf{1}; H) \xrightarrow{f_.} K(T_1; H) \to H/T_1H \to
0.
\end{displaymath}
Here $f_1=id$, $f_0=T_1$, and $H/T_1H$ is regarded as a complex
such that the $0$th component is $H/T_1H$, and other components
are zero.

Since $K(z_2, \cdots, z_{d})$ is a complex of free modules,
tensoring with it preserves exactness. Thus we obtain the
following exact sequence
\begin{displaymath}
0 \to K(1, z_2, \cdots, z_d; H) \to K(T_1, \cdots, T_d; H) \to
K(z_2, \cdots, z_d; H/T_1H) \to 0.
\end{displaymath}
But the tuple $(I, T_2, \cdots, T_d)$ is invertible, hence has
index zero. Now the conclusion follows from the additivity of
Euler characteristics. $\quad \square$

\bigno \textbf{Remark:} If the condition $ker(T_1)=\{0\}$ is
dropped, then Proposition \ref{P:reduction-index} is not true
which  can be seen by considering  $(M_z^*, M_w^*)$, the adjoints
of the multiplication operators on the Hardy space over the
bidisk.

\bigskip

Now we take up the Samuel multiplicity of  $H/T_1H$. First we
recall some general facts on Hilbert polynomials
(\cite{Hartshorne} page 49). While not strictly required, they
  provide the \emph{``the right way" } to present the material.
The reason largely lies in the identity
\begin{equation}\label{E:combin}
\left(
\begin{array}{c} z \\ d \end{array} \right) - \left(
\begin{array}{c} z-1 \\ d \end{array} \right) = \left(
\begin{array}{c} z-1 \\ d-1 \end{array} \right).
\end{equation}
Here $\left(
\begin{array}{c} z \\ r \end{array} \right) = \frac{1}{r!}\cdot
z(z-1)\cdots(z-r+1)$ is the binomial coefficient function.

 Note that Hilbert
polynomials are numerical polynomials. By a numerical polynomial
we mean a polynomial $P(z) \in \qq[z]$ such that $P(k) \in \zz$
for $k>>0$. Then for any numerical polynomial $P(z)$ of degree
$d$, there are integers $c_0, \cdots, c_d$, such that
\begin{equation}\label{E:expan}
P(z)=c_d \left(
\begin{array}{c} z \\ d \end{array} \right) +c_{d-1}\left(
\begin{array}{c} z \\ d-1 \end{array} \right) + \cdots + c_0.
\end{equation}
In particular, $P(k)\in \zz$ for all $k$. We call $c_d$ the
reduced leading coefficient of $P(z)$; it is   actually the
leading coefficient multiplied by $d!$. If we take
$P(k)=\phi_H(k)$ $(k
>>0)$ to be the Hilbert polynomial of a Hilbert module $H$, then
$c_d$ is just the Samuel multiplicity of $H$.

 For a function  $f$, defined on $\zz$, we define an operation
$\Delta f$ by taking the difference : $(\Delta f)(k)=f(k+1)-f(k)$.
Because of Equations (\ref{E:combin}) and (\ref{E:expan}), Hilbert
polynomials behave nicely under $\Delta$
\begin{equation*}
(\Delta P)(z)=c_d \left(
\begin{array}{c} z \\ d-1 \end{array} \right) +c_{d-1}\left(
\begin{array}{c} z \\ d-2 \end{array} \right) + \cdots + c_1.
\end{equation*}
 In particular, the
reduced leading coefficient is invariant under $\Delta$.

\bigskip

Fix a  Hilbert module $H$ over $A$, such that $dim \ H/IH <
\infty$, here $I=(z_1, \cdots, z_d)\subset A$. Let   $J=(z_2,
\cdots, z_d)\subset A$. Then $H/z_1H$ can be regarded as a module
over the polynomial ring in $d-1$ variables, say $z_2, \cdots,
z_d$, such that $dim \ (H/z_1H)/ J(H/z_1H) < \infty$. Similarly,
we define its Hilbert function
\begin{equation}
\phi_{H/z_1H}(k)= dim \ \frac{H/z_1H}{J^k(H/z_1H)}=dim \
\frac{H}{z_1H+J^kH},
\end{equation}
 and the Samuel multiplicity as the reduced
leading coefficient of $\phi_{H/z_1H}(k)$ at the degree $(d-1)$
term.

 Because of
$I^{k+1}H + z_1H=J^{k+1}H + z_1H$, the following is exact
\begin{equation*}
H/I^kH \xrightarrow{\tilde{z}_1} H/I^{k+1}H \xrightarrow{q}
H/(z_1H+J^{k+1}H) \to 0.
\end{equation*}
Here $\tilde{z}_1$ is induced by $z_1$, and $q$ is the natural
quotient map. It follows
\begin{proposition}\label{P:reduction-mul}
For any Hilbert module $H$ over the polynomial ring $\cc[z_1,
\cdots, z_d]$, such that $dim \ H/ IH < \infty$ here $I=(z_1,
\cdots, z_d)$, one has
\begin{equation*}
\phi_{H/z_1H}(k) \ge (\Delta \phi_H)(k)
\end{equation*}
 for all $k \in \nn$.
In particular,
\begin{equation*}
dim(H/IH) \ge e(J, H/z_1H) \ge e(I, H).
\end{equation*}
\end{proposition}

\subsection{More estimates of Hilbert
functions}\label{subS:estimate} It may appear natural to consider
the two variable homology Hilbert functions
\begin{equation}
h_i(s, t)=dim \ H_iK(T_1^s, T_2^t; H),
\end{equation}
for $i=0, 1, 2$, which, of  course, are not necessarily of
polynomial type, when $s, t >>0$. But even when $h_i(k)$, defined
by Equation (\ref{E:hik}), is a polynomial, $h_i(s, t)$ may still
fail to be so. It is largely open how to characterize when
$h_i(k)$ and $h_i(s, t)$ will be of polynomial type, when $k$,
$s$, $t$ are large enough. See \cite{TuCuong} for algebraic
results along this line.

 \bigno \textbf{Example} Let
$H^2(\mathbb{D})$ and $H^2(\mathbb{D}^2)$ be the Hardy space over
the unit disc and bidisc, respectively. Let $M_z, M_w$ be
multiplication by coordinate functions on $H^2(\mathbb{D}^2)$, and
$M_{\xi}$ on $H^2(\mathbb{D})$. Let
\begin{displaymath}
H=H^2(\mathbb{D}^2) \oplus H^2(\mathbb{D}), \quad T_1=\left(
\begin{array}{cc} M_z & 0 \\ P_0 & M_{\xi} \end{array} \right), \quad  and
\ \ T_2=\left(
\begin{array}{cc} M_w & 0 \\ P_0 & M_{\xi} \end{array} \right).
\end{displaymath}
Here $P_0: H^2(\mathbb{D}^2) \to H^2(\mathbb{D})$ is the
projection onto the constants, that is $P_0(f)=f(0, 0)$. Then
\begin{displaymath}
T_1T_2=T_2T_1=\left(
\begin{array}{cc} M_zM_w & 0 \\ M_{\xi} P_0 & M_{\xi}^2 \end{array}
\right), \qquad and \quad T_1^s=\left(
\begin{array}{cc} M_z^s & 0 \\ M_{\xi}^{s-1} P_0 & M_{\xi}^s \end{array}
\right).
\end{displaymath}
 So
\begin{equation*}
dim \ H/ (T_1^sH + T_2^tH)= st+ min\{ s, t \},
\end{equation*}
which is not a polynomial, but collapses to one when $s=t$. $\quad
\square$

\bigskip

 The main result of this subsection is Theorem \ref{T:cool}, which
 gives some
estimates  of   Hilbert functions  in more than one variable.
There exists an algebraic version of Theorem \ref{T:cool}
 due to C. Lech \cite{Lech1}.
The proof relies on the so called \emph{``form ideal"} technique,
which goes back to Krull \cite{Krull}.  Northcott's book
\cite{Northcott} contains a detailed discussion of this technique.
However, our proof is modelled after Lech \cite{Lech1}. It is
worth mentioning that, this technique of associating a homogenous
ideal to a module, has great computational power in calculating
Hilbert functions, especially when combined with Gr\"{o}bner bases
(\cite{Eisenbud}, chapter 15).

\begin{theorem}\label{T:cool}
Let $H$ be a Hilbert module  over the polynomial ring $\cc[z_1,
\cdots, z_d]$, such that $dim \ H/\overline{IH} < \infty$, here
$I=(z_1, \cdots, z_d)$, let $e=e(I,H)$ be the Samuel multiplicity
of $H$, we have
\begin{displaymath}
n_1\cdots n_d \cdot e \le dim \ H / (\overline{z_1^{n_1}H+\cdots +
z_d^{n_d}H}) \le n_1\cdots n_d(e + \frac{C}{min_{i} \ n_i})
\end{displaymath}
for some constant $C$ and any $n_1, \cdots, n_d \in \nn$.
\end{theorem}
Thus we obtain the following stronger form of Lech's formula

\begin{corollary}\label{C:Lech}
If $H$ is a Hilbert module over $\cc[z_1, \cdots, z_d]$ such that
$dim \ H/\overline{IH} < \infty$, then
\begin{equation*}
e(I, H)=\lim_{(min \ n_i) \to \infty} \frac{dim \
H/\overline{z_1^{n_1}H+ \cdots + z_d^{n_d}H}}{n_1 \cdots n_d}
\end{equation*}
\end{corollary}
If we replace $n_i$ by $n_i \cdot t$, and let $t \to \infty$, then
Part (iii) of Theorem \ref{T:main} follows.

\bigskip

Now specialize to $d=2$. By considering the adjoints $(T_1^*,
T_2^*)$, and Theorem \ref{T:cool}, we know that $h_i(s, t) \ge
s\cdot t \cdot e_i$ for $i=0, 2$. Now by Corollary \ref{C:dyna},
we have the following corollary, which will be used in the next
section.
\begin{corollary}\label{C:h1}
For any  pair $T=(T_1, T_2)$ of commuting operators, if $T$ is
Fredholm, and let $h_1(s, t)=dim \ H_1K(T_1^s, T_2^t)$, for $s$,
$t \in \nn$, then
\begin{equation*}
h_1(s, t) \ge s\cdot t \cdot e_1.
\end{equation*}
\end{corollary}
Recall that $e_1$ is defined in Theorem \ref{T:main}.

\bigno \textbf{Proof of Theorem \ref{T:cool}.}

Let $l=dim \ H/\overline{IH}$ (assuming $l>0$), and choose a
sequence of submodules
\begin{displaymath}
H=H_0\supset H_1 \supset \cdots \supset H_l=\overline{IH},
\end{displaymath}
such that $dim \ H_i/H_{i+1} =1$. This is always possible by
considering Jordan decomposition in linear algebra. Pick any
$\zeta_i \in H_i \setminus H_{i+1}$, such that $1 \to \zeta_i +
H_{i+1}$ induces an isomorphism $\cc \to H_i/H_{i+1}$.


Now multiply the above sequence by $I^\mu$, $\mu =0, 1, \cdots
k-1$, take closures in the Hilbert space $H$, and link them
together
\begin{displaymath}
H=H_0\supset H_1 \supset \cdots \supset H_l \supset
\overline{IH_1} \supset \overline{IH_2}\supset \cdots \supset
\overline{IH_l} \supset \cdots \supset \overline{I^{k-1}H_l}
(=\overline{I^kH}).
\end{displaymath}

Then for any ideal $\alpha \subset \cc[z_1, \cdots z_d]$, we can
add $\alpha H$ to the above sequence to obtain a chain from
$H=H+\alpha H$ to $\alpha H + \overline{I^kH}$. Note that $\alpha
H + \overline{I^kH}$ is automatically closed. Hence
\begin{equation}
dim \ H/(\alpha H + \overline{I^k H}) = \sum_{\mu=0}^{k-1}
\sum_{\nu=0}^{l-1} dim \ \frac{\alpha H
+\overline{I^{\mu}H_{\nu}}}{\alpha H + \overline{I^{\mu}H_{\nu +
1}}}.
\end{equation}

\bigno \emph{Fact:} For any vector spaces $A, M \supset N$, there
is a natural isomorphism between vector spaces
$\frac{A+M}{A+N}\cong \frac{M}{A\cap M+N}$.

\smallskip

\noindent \emph{Proof  of the Fact.} Note that $\frac{A+M}{A+N}
=\frac{M+A+N}{A+N}\cong \frac{M}{M\cap(A+N)}$. Then verify
directly that $M\cap(A+N) = M\cap A +N$.

\bigskip

Now we continue with the proof of Theorem \ref{T:cool}. It follows
that, as vector spaces
\begin{displaymath}
\frac{\alpha H +\overline{I^{\mu}H_{\nu}}}{\alpha H +
\overline{I^{\mu}H_{\nu + 1}}} \cong \frac{
\overline{I^{\mu}H_{\nu}}}{\alpha H \cap
\overline{I^{\mu}H_{\nu}}+ \overline{I^{\mu}H_{\nu + 1}}}.
\end{displaymath}

Let $F_{\mu}$ denote the collection of homogeneous polynomials of
degree $\mu$ in $\cc[z_1, \cdots, z_d]$. For any $\mu$ we define a
$\cc$-linear map $T_{\nu, \mu} : F_{\mu} \to
\frac{\overline{I^{\mu}H_{\nu}}}{ \overline{I^{\mu}H_{\nu + 1}}}$
by
\begin{equation*}
T_{\nu, \mu}(f)=f\zeta_{\nu} + \overline{I^{\mu}H_{\nu + 1}}.
\end{equation*}
Note that $z_j \zeta_{\nu} \in IH \subset H_{\nu + 1}$, and one
can then directly verify that, $T_{\nu, \mu}$ is surjective.

Moreover, it induces, through the quotient map, the following
\begin{equation*}
F_{\mu} \to \frac{\overline{I^{\mu}H_{\nu}}}{\alpha H \cap
\overline{I^{\mu}H_{\nu}} + \overline{I^{\mu}H_{\nu + 1}}},
\end{equation*}
whose kernel is denoted by
\begin{displaymath}
K_{\nu, \mu}(\alpha)=\{ f \in F_{\mu}, \ T_{\nu, \mu}f \subset
(\alpha H \cap \overline{I^{\mu}H_{\nu}}) +
\overline{I^{\mu}H_{\nu + 1}} \}
\end{displaymath}
Observe that $z_i K_{\nu, \mu}(\alpha) \subset  K_{\nu,
\mu+1}(\alpha)$, it follows that $I_{\nu}(\alpha)= \cup_{\mu
=0}^{\infty} K_{\nu, \mu}(\alpha)$ is a graded ideal in $\cc[z_1,
\cdots, z_d]$, called a \emph{``form ideal"} of $\alpha$. Then we
have
\begin{equation}
dim \ H/ (\alpha H + \overline{I^k H})=\sum_{\nu=0}^{l-1}
\sum_{\mu=0}^{k-1} dim \ \frac{F_{\mu}}{I_{\nu}(\alpha) \cap
F_{\mu}}.
\end{equation}
Two special cases are particularly important:

1. Let $\alpha=0$, then the number of $\nu$ such that
$I_{\nu}(0)=0$ is exactly $e(I, H)$, because, for a nonzero ideal
$I_{\nu}(0)$, the function $\sum_{\mu=0}^{k-1} dim \
\frac{F_{\mu}}{I_{\nu}(\alpha) \cap F_{\mu}}$ is a polynomial of
degree at most $d-1$ for large $k$;

 2. Assume that $\alpha$ is $I$-primary, that is, $I \supset \alpha \supset
I^t$ for some $t$. Then for  $k>t$,   one has
\begin{displaymath}
\alpha H + \overline{I^k H}=\overline{\alpha H +  I^k H}=
\overline{\alpha H}.
\end{displaymath}
 So
\begin{equation}\label{E:primary}
dim \ H / \overline{\alpha H} = \sum_{\nu=0}^{l-1} dim \
\frac{\cc[z_1, \cdots, z_d]}{I_{\nu}(\alpha)}.
\end{equation}

\bigskip

Now apply formula (\ref{E:primary}) to $\alpha=(z_1^{n_1}, \cdots,
z_d^{n_d})$, and observe that $\alpha \subset I_{\nu}(\alpha)$,
one has
\begin{equation*}
dim \ \frac{\cc[z_1, \cdots, z_d]}{I_{\nu}(\alpha)} \le n_1\cdots
n_d.
\end{equation*}

In addition, if $0 \ne f \in I_{\nu}(\alpha)$ for some $\nu$,
then, order the polynomials lexicographically, and assume that the
highest power in $f$ is $z_1^{\sigma_1}\cdots z_d^{\sigma_d}$.
Then, representatives of elements in $\frac{\cc[z_1, \cdots,
z_d]}{I_{\nu}(\alpha)}$ can be chosen in such a way that, any
power product $z_1^{r_1}\cdots z_d^{r_d}$ appearing in any
representative is subject to $r_i < n_i$, $i=1, 2, \cdots d$, and
in addition at least $r_i < \sigma_i$ for some $i$. Hence
\begin{equation*}
dim \ \frac{\cc[z_1, \cdots, z_d]}{I_{\nu}(\alpha)} \le n_1\cdots
n_d (\frac{\sigma_1}{n_1} + \cdots + \frac{\sigma_d}{n_d}).
\end{equation*}

Putting all $I_{\nu}(\alpha)$ together, we have an upper estimate
\begin{equation}\label{E:upper}
dim \ H / (\overline{z_1^{n_1}H + \cdots +z_d^{n_d}H}) \le
n_1\cdots n_d(e(I, H)+ \frac{A}{min_{i} \ n_i}),
\end{equation}
where $A$ is a constant independent of all $n_i$.

\bigskip

For the other direction of Theorem \ref{T:cool}, observe that, by
inequality (\ref{E:gee})
\begin{equation*}
dim \ H/(z_1^{n_1}H+\cdots + z_d^{n_d}H)  \ge   e((z_1^{n_1},
\cdots, z_d^{n_d}), H).
\end{equation*}
 So it suffices to show
\begin{equation*}
 e((z_1^{n_1},
\cdots, z_d^{n_d}), H) \ge n_1\cdots n_d \ e(I, H).
\end{equation*}

Note that the upper estimate (\ref{E:upper}) implies
\begin{equation}\label{E:in}
e(I, H) \ge \overline{\lim}_{(min \ n_i) \to \infty} \frac{dim \ H
/ (z_1^{n_1}H + \cdots +z_d^{n_d}H)}{n_1\cdots n_d}.
\end{equation}
Choose $m_1, \cdots, m_d$ such that $n_1m_1=\cdots=n_dm_d$, then
apply the above inequality (\ref{E:in}) to the ideal $(z_1^{n_1},
\cdots, z_d^{n_d})$, instead of $I$ to obtain
\begin{eqnarray*}
  e((z_1^{n_1},
\cdots, z_d^{n_d}), H)  & \ge &  \overline{\lim}_{k\to \infty}
\frac{dim \ H/((z_1^{n_1})^{m_1k}H+\cdots +
(z_d^{n_d})^{m_dk}H)}{(m_1k)\cdots(m_dk)} \\
& \ge & n_1\cdots n_d \lim_{t\to \infty} \frac{dim \
H/(z_1^tH+\cdots + z_d^tH)}{t^d}.
\end{eqnarray*}
By  Theorem \ref{T:H-Lech} (Lech's formula), the last expression
is just $n_1\cdots n_d \ e(I, H)$. Now our proof of Theorem
\ref{T:cool} is complete. $\quad \square$

\section{Fredholm theory on the symmetric Fock space in two variables}
In this section we show that  the algebraic theory developed in
Section 1 provides useful tools  in the study of operator theory
and function theory on the symmetric Fock space over the ball in
$\cc^2$. In particular, our results imply that Arveson's program
\cite{Arveson-problem} on Fredholm theory for pure d-contractions
can have satisfactory answers for both the Fredholmness problem
and the calculation of Fredholm indices when $d=2$. It is
noteworthy to mention that the Fredholmness problem is a subtle
one in general.

\medskip

In multivariable Fredholm theory the lack of amenable examples
  has seriously hampered the development of the subject. When trying to obtain
   meaningful examples, it is a folk problem to determine the
  Fredholmness
   and the index of the tuples of multiplication operators on holomorphic function spaces.
For example, a complete characterization of the Fredholmness of
submodules of the vector-valued Hardy module over the polydisc
$H^2(\dd^n) \otimes \cc^N$, and to calculate the index, seem out
of reach at this time. This might require the development of more
sophisticated machineries from both the analytic side, that is,
function theory on the polydisc, and the algebraic side, that is,
theory of Koszul complexes over Hilbert modules. When $n=2$,
$N=1$, Yang \cite{Yang} showed that the index of the pair of
multiplication by coordinate functions is one for a fairly large
class. See Curto-Salinas \cite{Curto-Salinas} for a study of
Bergman-type spaces. Also see \cite{E-Putinar} for more
discussions on Bergman spaces. Through Arveson's work
\cite{subalgebra}, it is now known that Fredholm theory on the
symmetric Fock space is in fact equivalent to that of pure
d-contractions. In particular, in a series of paper \cite{Dirac},
\cite{Arveson-psum}, \cite{Arveson-problem}, Arveson initiated a
study of Fredholm theory for d-contractions. In \cite{Gleason}
Gleason-Richter-Sundberg  showed that, under certain regularity
condition on the invariant subspace of the vector valued symmetric
Fock space, the associated Koszul complex (in possibly more than
two variables) is acyclic, that is, all homology groups are zero
except for the last stage, hence the index is equal to the
dimension of the last homology group.

\medskip

In general, the difficulty in calculating the several variable
Fredholm index is, mainly due to the lack of an effective way to
estimate the intermediate homology groups, that is, those other
than the first and last. A common phenomenon, which has occured in
calculating the index in past research,  is that the Koszul
complex $K(\cdot)$ is shown to be acyclic, often by nontrivial
arguments, and hence the index is identified with $dim \
H_0K(\cdot)$.

\medskip

 Now recall that the symmetric
Fock space over the ball $B_2 \subset \cc^2$, denoted by $H^2$, is
a Hilbert space of holomorphic functions, defined on the ball
$B_2\subset \cc^2$. It is determined by the reproducing kernel
$k((z,w), (\zeta, \eta))=\frac{1}{1-z\bar{\zeta}-w\bar{\eta} }$;
equivalently, it can be obtained by symmetrizing the full Fock
space, as was done by Arveson \cite{subalgebra}. In what follows,
let $\submo \subset H^2 \otimes L$ be a multiplication invariant
subspace of the $L$-valued symmetric Fock space, here $L$ is a
Hilbert space; and $(M_z, M_w)$ be the pair obtained by
restricting   multiplication by coordinate functions $z, w$ to
$\submo$.
Note that $H_2K(M_z, M_w; \submo)=0$, and $H_0K(M_z, M_w;
\submo)=\submo/(z\submo +w\submo)$. So, in order to have a
Fredholm pair, an obvious necessary condition is to have $dim \
\submo/(z\submo +w\submo)< \infty$, which we shall show is also
sufficient.

The following is the main result of this section.
%
\begin{theorem}\label{T:subFock}
Let $\submo \subset H^2\otimes L$ be a multiplication invariant
subspace of the $L$-valued symmetric Fock space in two variables.
Here $L$ is a Hilbert space. Let $(M_z, M_w)$ be the pair of
multiplication by coordinate functions  on $\submo$.
\begin{itemize}
\item [(1)] $(M_z, M_w)$ is Fredholm if and only if $dim \
\submo/(z\submo + w\submo) < \infty$.

\item [(2)] When Fredholm, the index is
\begin{equation*}
index(M_z, M_w)=e_0 \ (=Samuel \ multiplicity \ of \ \submo).
\end{equation*}
\end{itemize}
\end{theorem}
\textbf{Proof.}   We shall try to separate the algebraic and
analytic components in our proof as much as possible. This is done
through step 1 - 7 below.   Basically there are three main
ingredients: 1. general algebraic properties of Koszul complexes;
2. Fredholm theory developed in Section 1 and \cite{Fred}; and 3.
function-theoretic operator theory, such as reproducing kernels,
holomorphic multipliers, and the Hardy space over the unit disc
$H^2(\dd)$.

\medskip

We first set two notations. For any set $S$, consisting of
functions in two variables, let $S(a, b)=\{f(a, b), f\in S\}$. For
any holomorphic function $f$ in one variable, defined on a
neighborhood of the origin, if $f(x)=c_rx^r +c_{r+1}x^{r+1}+
\cdots$ is its Taylor series, and $c_r\ne 0$, then we define
$ord_x \ f =r$.

\bigno \emph{Step 1.}  Let $T_w$ be the map acting on the vector
space $\submo/z\submo$, induced by   multiplication by $w$. Then
the two variable index is reduced to the one variable index
\begin{lemma}\label{L:2to1}
For any invariant subspace $\submo$, one has
\begin{equation}
index(M_z, M_w; \submo)=-index(T_w; \submo/z\submo).
\end{equation}
In particular, $(M_z, M_w)$ is Fredholm ( semi-Fredholm, not
Fredholm, resp.) if and only if $T_w$ is Fredholm ( semi-Fredholm,
not Fredholm, resp.).
\end{lemma}
Here the latter index is defined in the following more general
setting: Let $T$ be a linear map on a vector space $V$, then   its
index is $index(T; V)=dim \ ker(T) - dim \ V/TV$, if finite.

\bigno \textbf{Proof.} This result follows from Proposition
\ref{P:reduction-index}. Since we shall need more precise
information between $(M_z, M_w)$ and $T_w$   in what follows, we
give a different approach. For any pair $(T_1, T_2)$ acting on
$H$, we consider the following exact sequence of complexes
\begin{displaymath}
0 \to K(T_1; H) \to K(T_1, T_2; H) \to \bar{K}(T_1; H) \to 0,
\end{displaymath}
 where $\bar{K}$ is the complex obtained by shifting $K$ to
the right by one step. We have the following long exact sequence
\begin{displaymath}
0 \to H_2K(T_1, T_2; H) \to H_1K(T_1; H)
\xrightarrow{-\tilde{T}_2} H_1K(T_1; H) \to H_1K(T_1, T_2; H)
\end{displaymath}
\begin{displaymath}
\to H_0K(T_1; H) \xrightarrow{-\tilde{T}_2}  H_0K(T_1; H) \to
H_0K(T_1, T_2; H)   \to 0.
\end{displaymath}
Here $\tilde{T}_2$ is the map induced by $T_2$. When applied to
$(M_z, M_w)$ on $\submo$, the above sequence yields
\begin{equation}
0 \to H_1K(M_z, M_w; \submo) \to \submo/z\submo \xrightarrow{T_w}
\submo/z\submo \to H_0K(M_z, M_w; \submo) \to 0.
\end{equation}

In particular, we know that $ker(T_w)\cong H_1K(M_z, M_w;
\submo)$. Similarly, it follows that $ker(T_w^k)\cong H_1K(M_z,
M_w^k; \submo)$, for any $k \in \nn$, which will be used in step
7.

\bigno \emph{Step 2.}
\begin{lemma}
For any pair $(T_1, T_2)$ of commuting operators on a Hilbert
space $H$, if $U=\left( \begin{array}{cc} a & c
\\ b & d \end{array} \right)$ is a unitary matrix, then
for $i=0, 1, 2$,
\begin{equation}
H_iK(T_1, T_2; H) \cong H_iK(aT_1+bT_2, cT_1+dT_2; H)
\end{equation}
are natural  isomorphisms between vector spaces.
\end{lemma}
The proof is essentially linear algebra, and is best carried out
in the \emph{``exterior power definition"} of Koszul complexes.
Readers can show that the matrix $U$ actually induces an
isomorphism between the Koszul complexes of $(T_1, T_2)$ and
$(aT_1+bT_2, cT_1+dT_2)$. The details are skipped.

\bigskip

So in order to consider the Fredholm theory of $(T_1, T_2)$, it is
equivalent to look at  $(aT_1+bT_2, cT_1+dT_2)$. This is
particularly applicable to \emph{``ball-oriented operator
theory"}. Note that $(z, w) \to (z', w')=(az+bw, cz+dw)$ is a
change of variables on the ball $B_2$, which preserves the metric.
Moreover, it preserves the metric on $H^2$, that is, for any two
variable polynomial $p(\cdot, \cdot)$, one has
\begin{displaymath}
||p(z, w)||_{H^2}=||p(az+bw, cz+dw)||_{H^2},
\end{displaymath}
which can be seen easily from the reproducing kernel: the change
of variables $(z, w) \to (z', w')=(az+bw, cz+dw)$ does not change
the form of the kernel $\frac{1}{1-z\bar{\zeta}-w\bar{\eta} }$,
because $\left(
\begin{array}{cc} a & c
\\ b & d \end{array} \right)$ is unitary, hence if
$(\zeta', \eta') = (a\zeta + b\eta, c\zeta+ d\eta)$, then
$z\bar{\zeta}+w\bar{\eta}= z'\bar{\zeta'}+w'\bar{\eta'}$.

 So the pair of multiplication operator $(M_{az+bw}, M_{cz+dw})$ on
an invariant subspace $\submo \subset H^2\otimes L$ is unitarily
equivalent to $(M_z, M_w)$ on another invariant subspace $\submo'
\subset H^2\otimes L$, obtained from $\submo$ by a change of
variables.


\bigno \emph{Step 3.} Now we need the Fredholm theory developed in
\cite{Fred}. We rephrase Theorem 2 in \cite{Fred} in the setting
of a general linear map $T: V \to V$, acting on a vector space.
Define the shift and backward shift Samuel multiplicity of $T$ by
\begin{equation*}
\smult =\lim_{k \to \infty}\frac { dim \ V/T^kV}{k}, \quad  and \
\ \bsmult=\lim_{k \to \infty}\frac { dim \ ker(T^k)}{k}.
\end{equation*}
 If $T$ is
semi-Fredholm in the general sense, that is, if at least one of
$dim \ ker(T)$ and $dim \ V/TV$ is finite, then
\begin{itemize}
\item[(i)]$\smult, \ \bsmult \in \{0, 1, 2, \cdots, \infty\}$;

\item[(ii)] $index(T)= \bsmult - \smult$;

\item[(iii)]when $k>>0$,
\begin{equation*}
\bsmult=dim(ker(T)\cap T^kH)=dim(ker(T) \cap T^{\infty}H),
\end{equation*}
and
\begin{equation*}
\smult=dim(\frac{H}{TH+ker(T^k)})=dim(\frac{H}{TH+ker(T^{\infty})}).
\end{equation*}
\end{itemize}
Here $T^{\infty}H=\cap_{k \ge 0} T^kH$, and
$ker(T^{\infty})=\cup_{k \ge 0} ker(T^k)$.

\bigskip

The proof in \cite{Fred} is for Hilbert space operators, but a
careful examination of the proof will reveal that  it carries over
to the more general setting of vector spaces. Only a small change
is needed: In \cite{Fred}, $ker(T^{\infty})$ is defined to be
$\overline{\cup_{k \ge 0} ker(T^k)}$. Here, we do not have a
Hilbert space structure, and hence cannot form the closure. On the
other hand, $ker(T^{\infty})$ appears only in
$TH+ker(T^{\infty})$.
By arguments in \cite{Fred}, it is easy to see that $TH+ker(T^k)=
TH+ker(T^{\infty})$, for large $k$, and is always closed for a
semi-Fredholm map in the general sense. So, it does not matter if
we do not form the Hilbertian closure of $ker(T^{\infty})$. Recall
that $ker(T) \cap T^{\infty}H$ is defined to be  the
\textbf{stabilized kernel} of $T$ \cite{Fred}.

\bigno \emph{Step 4. }  Now let $V=\submo/z\submo$, and we apply
step 3 to the linear map $T_w: V \to V$ to calculate its Fredholm
index by looking at the stabilized kernel and cokernel.

According to McCullough-Trent \cite{McCullough}, we can assume
that, there is a partial isometric multiplier $\Phi : H^2 \otimes
E \to H^2 \otimes L$, where $E$ is a Hilbert space, such that
$\Phi(H^2 \otimes E)=\submo$. Let $\submo'=ker(\Phi)$.

Observe that
\begin{displaymath}
T_w^{\infty}V\cong \frac{\cap_{k \ge 0}(w^k \submo
+z\submo)}{z\submo}.
\end{displaymath}

\begin{lemma}
For any $x=\Phi \xi \in \submo$, we have
\begin{itemize}
\item[(i)] $x \in \cap_{k \ge 0}(w^k \submo +z\submo)$ if and only
if
\begin{equation*}
 \xi(0, w) \in \cap_{k \ge 0} (w^k H^2(0, w) \otimes E
+\submo'(0, w)).
\end{equation*}


\item[(ii)] $x \in z \submo$ if and only if $\xi(0, w) \in
\submo'(0, w)$.
\end{itemize}
\end{lemma}
Proof. Part (i) : The following statements  are all equivalent.

1. $\Phi \xi \in \cap_{k \ge 0}(\Phi (w^k H^2\otimes E +z H^2
\otimes E))$;

2. for any $k$, there exists $u, v \in H^2\otimes E$, such that
$\Phi(\xi -w^ku -zv)=0$;

3.  $\xi \in w^kH^2\otimes E+zH^2 \otimes E +\submo'$, for any
$k$;

4.  $\xi(0, w) \in w^k H^2(0, w) \otimes E +\submo'(0, w)$, for
any $k$.

\bigskip

For Part (ii), $\Phi \xi \in z\submo=z\Phi (H^2 \otimes E)$ means
$\Phi(\xi-z\eta)=0$, for some $\eta \in H^2 \otimes E$. That is,
$\xi \in \submo' +zH^2 \otimes E$, which is equivalent to $\xi(0,
w)\in \submo(0, w)$. $\quad \square$

\bigno \emph{Step 5.}  Next we deal with $\submo'(0, w)$, which is
a subspace of $H^2(0, w)\otimes E$, which in turn is equivalent to
the direct sum of $dim(E)$ many copies of the Hardy space
$H^2(\dd) $ over the unit disc $\dd \subset \cc$. Now the
difficulty is that, we do not know whether $\submo'(0, w)$ is
closed.

Observe that $\Phi(0, w) H^2(0, w)\otimes E =\submo'(0, w)$.
Thanks to the result of Greene-Richter-Sundberg \cite{Greene}, we
know that there exists a (thin) subset $Z \subset S=\partial B_2
\subset \cc^2$ of the sphere, which is contained in the zero set
of a nonzero holomorphic function, such that the boundary values
of the partial isometric holomorphic multiplier $\Phi$ are partial
isometries,
  with a constant rank on $S \setminus Z$.  Now by arguments in
  step 2,
   we can apply a change of variables, if necessary, such
  that $Z \cap \{(0, w), \ |w|=1\}$ is a thin subset of the circle
$\{(0, w), \ |w|=1\}$. It follows that $\Phi(0, w)$ is a
one-variable, bounded, holomorphic multiplier, with boundary
values that are partial isometries with   constant rank almost
everywhere on the unit circle.

Now by the familiar theory of $H^2(\dd)$, we know that $\Phi(0,
w)$, acting on the direct sum of Hardy spaces $H^2(0, w)\otimes
E$, is a partial isometry, and hence has closed range. That is, we
can assume that  $\submo'(0, w)$ is closed.

\bigno \emph{Step 6.}  In this step we show that the stabilized
kernel of $T_w$ is $0$. That is to show that
\begin{quote}
if $x \in \cap_{k \ge 0} (w^k\submo +z\submo)$, and $wx \in
z\submo$, then $x\in z\submo$.
\end{quote}
 Let $x =\Phi \xi$, where $\Phi$ is from step 4. Then by   step 4 we need to show that
 \begin{quote}
if $ \xi(0, w) \in \cap_{k \ge 0} (w^k H^2(0, w) \otimes E
+\submo'(0, w))$, and $ w\xi(0, w) \in \submo'(0, w)$, then $
\xi(0, w) \in \submo'(0, w)$.
\end{quote}
%
By the Beurling-Lax-Halmos Theorem   on the Hardy space
$H^2(\dd)$, there exists an isometric (not just partially
isometric) holomorphic multiplier $\Theta$, in terms of $w$, such
that $\Theta (H^2(\dd) \otimes F) =\submo'(0, w)$ for some Hilbert
space $F$. For each $k \ge 0$ we can write $\xi(0, w)=w^k h_k
+\Theta g_k$ for some $h_k \in H^2(\dd) \otimes E$, and $ g_k \in
H^2(\dd) \otimes F$. Moreover, we write $w\xi(0, w) = \Theta g$,
where $g \in H^2(\mathbb{D}) \otimes F$.

We first show that the sequence $h_k$ can be chosen in such a way
that $sup_{k \ge 0} \ ||h_k|| < \infty$. To achieve this, we write
$g=g(0) + w g'$, and let $T_k(g')$ denote the Taylor polynomial of
$g'$ at the origin with degree $k$. Now, if we let $g_k=T_k(g')$,
then by using the above expressions of $\xi(0, w)$ and $w\xi(0,
w)$, we  have
\begin{equation}\label{E:fix}
w^{k+1}h_k= \Theta g(0) +w \Theta (g'-T_k(g')).
\end{equation}
 It follows that
for this choice of $h_k$ one has  $\lim_{k\to \infty}
||h_k||=||\Theta g(0)||$.

Now we form the inner product between $\Theta g(0)$ and both sides
of Equation (\ref{E:fix}), and observe that $\Theta g(0)$ and $w
\Theta (g'-T_k(g'))$ are orthogonal, we obtain
\begin{equation*}
||\Theta g(0)||^2 = \langle \Theta g(0), w^{k+1}h_k \rangle.
\end{equation*}
Since $||h_k||$ is uniformly bounded, the right hand side tends to
zero as $k \to \infty$. If follows $\Theta g(0) =0$. Recall that
$\Theta$ is an isometry, so $g(0)=0$. Now $\xi(0, w) = \Theta g'
\in \submo'$.

 \bigno \textbf{Remark:} The above arguments fail  for
the Bergman space: Take an invariant subspace $\submo \subset
L^2_a(\dd)$ of the Bergman space over the unit disc, such that
$dim \ \submo / w \submo = 2$ and $\submo(0)\ne \{0\}$. Then there
exists a function, say h, in $\submo \ominus w \submo$ such that
$h(0)=0$. So if we write $h=wh' \in \submo$, then it does not
follow that $h' \in \submo$.

\bigno \emph{Step 7.} Now we are ready to complete the proof of
Theorem \ref{T:subFock}. Recall from \cite{Fred} that, the
difference between  the kernel $ker(T_w)$ and  the stabilized
kernel $ker(T_w)\cap T_w^{\infty}V$ is
\begin{displaymath}
\frac{ker(T_w)}{ker(T_w)\cap T_w^kV} \cong
\frac{ker(T_w^k)+T_wV}{T_wV}
\end{displaymath}
for large $k$, which is  finite dimensional for a semi-Fredholm
map. This is where we use the condition $dim \ \submo/(z\submo
+w\submo) < \infty$. Since the stabilized kernel of $T_w$ is $0$,
it follows that $ker(T_w)\cap T_w^kV =0$ for large $k$.  So we
know that $ker(T_w) \cong \frac{ker(T_w^k)+T_wV}{T_wV}$ is finite
dimensional, and hence $T_w$ is Fredholm.

 By results in \cite{Fred} or step 3, the stabilized kernel of $T_w$
is $0$ which  means that
\begin{equation*}
sup_{k \ge 0} \ dim \ ker(T_w^k) < \infty.
\end{equation*}

 By the remark at the end of the proof of Lemma \ref{L:2to1}, this means
\begin{equation*}
sup_{k \ge 0} \ h_1(1, k) < \infty.
\end{equation*}

 By Corollary \ref{C:h1}, we know
that $h_1(t, k)\ge t \cdot k \cdot  e_1$, for any $t, k \ge 0$. It
follows that $e_1=0$.

Obviously, $e_2=0$. So by Theorem \ref{T:main}, we have
$index(M_z, M_w)=e_0$, and   the proof of Theorem \ref{T:subFock}
is complete.


\section{Additivity of Hilbert polynomials}
\subsection{Additivity and monotonicity of Samuel
multiplicities}\label{subS:add} For any invariant subspace $\submo
\subset H^2(\mathbb{D})$ of the Hardy space one has the well known
codimension-one property $dim \ \submo \ominus z \submo =1$. This
result has been exploited in   one variable, but a multivariable
generalization proves to be  resistent.  At the end of
\cite{Dirichlet} it is conjectured that the right approach is
probably through considering certain stabilized dimensions. In
this section, as a consequence of some algebraic considerations on
Hilbert spaces, we prove what might be called the generalization
of the codimension-one property for the symmetric Fock space in
two variables (Corollary \ref{C:mono}).

\medskip

 As we have seen, the study of
multivariable Fredholm index is closely related to the study of
Hilbert polynomials. Hilbert polynomials play essential roles in
algebraic geometry. One of their key properties is their
additivity over short exact sequences, which can in fact
characterize the Hilbert polynomials in some sense (see
\cite{Eisenbud} Exercise 19.18).

For a short exact sequence of finitely generated modules over a
Noetherian ring $R$
\begin{displaymath}
0\to L\to M \to N \to 0,
\end{displaymath}
with $\phi_{L}$, $\phi_{M}$, and  $\phi_{N}$ the Hilbert
polynomials of $L, M$, and $N$ respectively with respect to an
ideal $I\subset R$, the additivity of Hilbert polynomials in the
graded case refers to
\begin{displaymath}
\phi_{L}+\phi_{N}=\phi_{M};
\end{displaymath}
In  \cite{Eisenbud} Exercise 19.18,  it is shown that the Hilbert
functions are the universal additive invariants of graded modules
over the polynomial ring $\cc[z_0, \cdots, z_d]$, or equivalently,
they are the universal additive invariants of coherent sheaves on
the complex projective space $\mathbb{P}^d$. In a recent paper
\cite{Chan}, Chan showed that  this additivity of Hilbert
polynomials is equivalent to the multiplicativity of Chern numbers
over   projective spaces.

In the non-graded case, the additivity of the  leading terms, that
is, the additivity of Samuel multiplicities,
 is still true under mild conditions. In \cite{Dirichlet} such an
 additivity is shown to be true on the Dirichlet space over the
 unit disc, which has several consequences in operator theory.

\bigskip

 The purpose of this subsection is to establish a
version of the additivity of Hilbert polynomials for $H^2$
(Theorem \ref{T:add}).  Theorem \ref{T:subFock} plays an important
role in the proof of Theorem \ref{T:add}. We conjecture that the
conclusion of Theorem \ref{T:add} still holds for three or more
variables. The problem is open on the Hardy space over the
polydisc or the ball in $\cc^d$ $(d \ge 2)$.
\begin{theorem}\label{T:add}
Let $\submo \subset H^2 \otimes \cc^N$ $(N\in \nn)$ be an
invariant subspace of the $\cc^N$-valued symmetric Fock space in
two variables; let $ \submo^{\perp}=H^2 \otimes \cc^N \ominus
\submo$. Equip   $\submo$, $H^2 \otimes \cc^N$, and
$\submo^{\perp}$ with the natural module structures over $\cc[z,
w]$. Let $I=(z, w)$.

If $dim \ \submo/I \submo < \infty$, then
\begin{equation*}
e(I, \submo)+e(I, \submo^{\perp})=N.
\end{equation*}
\end{theorem}

Observe that, in order to formulate the above additivity of
Hilbert polynomials, it is necessary to have the condition $dim \
\submo/I \submo < \infty$, so that the Samuel multiplicity $e(I,
\submo)= d! \ \lim_{k \to \infty} \frac{dim \ \submo/I^k
\submo}{k^d}$ is finite, where $d=2$. Also, we do not form the
closure $\overline{I \submo}$, in order to match up with the
algebraic formulation.

We have the following \emph{``monotonicity property"} of $e(I,
\submo)$.
\begin{corollary}\label{C:mono}
If $\submo_1 \subset \submo_2 \subset H^2 \otimes \cc^N$ $(N\in
\nn)$ are two invariant subspaces of the $\cc^N$-valued symmetric
Fock space in two variables, such that $dim \ \submo_1/I \submo_1
< \infty$, and  $dim \ \submo_2/I \submo_2 < \infty$, then
\begin{displaymath}
e(I, \submo_1) \le e(I, \submo_2).
\end{displaymath}
In particular, for any invariant subspace $\submo \subset H^2
\otimes \cc^N$ $(N\in \nn)$, such that $dim \ \submo/I \submo <
\infty$, one has
\begin{displaymath}
e(I, \submo) \le N.
\end{displaymath}
\end{corollary}
\textbf{Proof.} Let $F_n$ denote the polynomials in $H^2 \otimes
\cc^N $ of degree at most $n$. Observe that
\begin{displaymath}
dim \ \submoperp/I^n \submoperp =dim \ \cap_{\alpha_1+\alpha_2=n}
\ ker({M_z^{*}}^{\alpha_1}{M_w^{*}}^{\alpha_2})=dim \ \submoperp
\cap F_n,
\end{displaymath}
where $\submo \subset H^2 \otimes \cc^N$ is a general invariant
subspace. Since
\begin{equation*}
\submo_{1}^{\perp}\cap F_n \supset \submo_{2}^{\perp}\cap F_n,
\end{equation*}
the definition of  the Samuel multiplicity implies $e(I,
\submo_1^{\perp}) \ge e(I, \submo_1^{\perp})$. An application of
Theorem \ref{T:add} completes the proof. $\quad \square$

\bigno \textbf{Remark} The above corollary is somewhat perplexing,
in the sense that, for an invariant subspace $\submo \subset H^2
\otimes \cc^N$ with $dim \ \submo/I \submo < \infty$, its
sub-invariant subspace $\submo' \subset \submo$ may still have
$dim \ \submo'/I \submo' =\infty$. Hence $e(I, \submo')=  d! \
\lim_{k \to \infty} \frac{dim \ \submo'/I^k \submo'}{k^d}=\infty$,
where $d=2$. But when   its Samuel multiplicity  is finite, it is
dominated by that of $\submo$.

If we replace $H^2$ by the one variable Hardy or Dirichlet space
over the unit disc, then the second part of Corollary \ref{C:mono}
is the well known codimension-N-property. The straight forward
generalization of this property to higher dimensions, obtained by
considering the codimension $dim \ \submo/I\submo$, fails
immediately; for instance, if $\submo$ is the invariant subspace
$[z, w]$ of $H^2$, consisting of functions vanishing at the
origin, then $dim \ \submo/I\submo =2 > 1$. In Section 7 of
\cite{Dirichlet} an explanation of this phenomenon was given, and
it was suggested that a notion of   \emph{``stabilized
codimension"} might be what is really relevant. It is worth
mentioning that, in the one variable case, the codimension is
always stabilized. Corollary \ref{C:mono} provides such a
\emph{``stabilized codimension-N-property"} in more than one
variable.

 \bigno \textbf{Proof of Theorem \ref{T:add}. }  We first show
 that $e(I, \submoperp)$ is always well defined, even when $\submo
 \subset H^2 \otimes \cc^N$ is not Fredholm. Then we find the
 corresponding invariant on $\submo$, denoted by $\sigma=\sigma_{\submo}$,
which is additive with respect to $e(I, \submoperp)$ on
$\submoperp$, \emph{i.e.},  $\sigma+e(I, \submoperp)=N$.

 To achieve this
it is common
  in commutative or homological algebra to consider the
exactness of the functor sending the module $H=\submoperp$ to
$H/I^kH$, which is in fact right half-exact. Hence with respect to
the exact sequence of Hilbert modules $0\to \submo \to H^2 \otimes
\cc^N \to H \to 0$,  the following sequence is exact
\begin{equation}
\to \frac{\submo}{I^k \submo} \to \frac{H^2 \otimes \cc^N}{I^kH^2
\otimes \cc^N} \to \frac{H}{I^kH} \to 0.
\end{equation}
Of course,  the exactness of the above sequence can be verified
directly. Since the above middle term is finite dimensional, so is
$H/I^kH$. It follows that $I^kH$ is closed in $H$, and $e(I,
\submoperp)$ exists, and is finite.

To complete the above sequence, we look at the image of the second
arrow and it follows that the following sequence is exact
\begin{equation}\label{E:unknown}
0 \to \frac{(\submo+ I^kH^2 \otimes \cc^N)}{I^kH^2 \otimes \cc^N}
\to \frac{H^2 \otimes \cc^N}{I^kH^2 \otimes \cc^N} \to
\frac{H}{I^kH}\to 0.
\end{equation}
Observe that,
\begin{displaymath}
dim \ \frac{(\submo+ I^kH^2 \otimes \cc^N)}{I^kH^2 \otimes \cc^N}
=dim\ P_{k-1}\submo,
\end{displaymath}
where   $P_k$ is the projection onto the space of polynomials of
degree $\leq k$. By the above exact sequence (\ref{E:unknown}),
the function
\begin{equation}
\varphi_{\submo}(k)=dim \ P_{k-1}\submo
\end{equation}
 will be a polynomial of
degree at most $d$ for $k>>0$, and $d! \lim_{k\to
\infty}\frac{\varphi_{\submo}(k)}{k^d}$ is an integer. Here $d=2$.
For convenience, we define
\begin{equation}\label{E:sigma}
(\sigma=)\sigma_{\submo}=d! \lim_{k\to \infty}\frac{dim\
P_{k-1}\submo}{k^d}.
\end{equation}
 Then $\sigma + e(I, \submoperp)=N$.

 \bigskip

 Next we take up $e(I, \submo)$. From Theorem
 \ref{T:subFock},
 we know that $(M_z, M_w)$ is
Fredholm, and $e(I, \submo)$  is equal to $index(M_z, M_w)$. Now
we need a result of Gleason-Richter-Sundberg \cite{Gleason}: if
$(M_z, M_w)$ is assumed to be Fredholm, then its index is equal to
the fibre dimension $f.d.(\submo)$. Here the fibre dimension is
defined by
\begin{equation}\label{E:fibre}
f.d. (\submo) = sup \ \{ dim \ \submo(z, w), \ |z|^2+|w|^2 <1\}.
\end{equation}
Note that the sup is achieved almost everywhere.

Now, what we need to show is that
\begin{quote}
$ (\ast) \qquad \qquad $ if $(M_z, M_w)$ is Fredholm, then $f.d.
(\submo) = \sigma$.
\end{quote}
 Here
it is worthwhile to point out that, the definition of both $f.d.
(\submo)$ and $ \sigma$ (see  (\ref{E:sigma}), and
(\ref{E:fibre})), does not require  the Fredholmness condition.
Hence, statement $(\ast)$ might be true   more generally. This
leads us to   Theorem \ref{T:identify}, which is true for the
symmetric Fock space $H^2_d$ in $d$ variables, $d\ge 2$.
\begin{theorem}\label{T:identify}
Let $\submo \subset H^2_d \otimes \cc^N$ be an invariant subspace
of  the  $\cc^N$-valued symmetric Fock space in $d$ variables
$(d\in \nn)$. Then
\begin{equation*}
f.d.(\submo)= d! \lim_{k\to \infty}\frac{dim \ P_k {\submo}}{k^d}.
\end{equation*}
\end{theorem}
%
The proof of Theorem \ref{T:identify} is in Subsection
\ref{subS:identifyproof}, which will complete   the proof of
Theorem \ref{T:add}. A generalization of Theorem \ref{T:identify}
will be given in \cite{index-of-quotient}. Before giving the proof
we list several other consequences of Theorem \ref{T:identify} in
Subsection \ref{subS:curvature}, mostly having something to do
with the curvature invariant introduced by Arveson
\cite{curvature}.

\subsection{The curvature of a pure
d-contraction}\label{subS:curvature} This subsection contains
applications to Arveson's curvature invariant. This is largely in
response to Arveson's question on expressing the curvature
invariant by invariants which are directly determined by the
spatial actions of the d-contractions, and through the expression
one can immediately tell that the curvature is an integer
\cite{Dirac}, \cite{Arveson-psum}. The content of this and the
next subsections, together with Theorem \ref{T:identify}, has been
circulated in a preprint under the title \emph{``Samuel
multiplicity and Arveson's curvature invariant"}.

\medskip

 Recall that \cite{subalgebra}
a d-contraction is a d-tuple of commuting operators
$T=(T_1,\cdots, T_d)$ acting on a Hilbert space $H$ that defines a
row contraction in the sense that
\begin{displaymath}
||T_1\xi_1+\cdots+T_d\xi_d||^2 \leq ||\xi_1||^2+\cdots+||\xi_d||^2
\end{displaymath}
for all $\xi_1, \cdots, \xi_d \in H$. For every d-contraction we
have $T_1T_1^*+\cdots+T_dT_d^*\leq \textbf{1}$. Define the defect
rank of the d-contraction $T$ to be the rank of the defect
operator $\Delta_T=\sqrt{\textbf{1}-T_1T_1^*-\cdots-T_dT_d^*}$.
$T$ is said to be pure if the completely positive map defined by
\begin{equation*}
\psi(X)=T_1XT_1^*+\cdots+T_dXT_d^*, \quad X\in B(H)
\end{equation*}
satisfies $\psi^n(\textbf{1})\to 0$ strongly, as $n\to \infty$. In
the case of a single contraction, it reduces to the condition that
$T$ belongs to the class $C_{., 0}$, see Sz.-Nagy and Foias
\cite{Nagy}.

For any $z\in \cc^d$,   define
$T(z)=\bar{z}_1T_1+\cdots+\bar{z}_dT_d$, and
$F(z)=\Delta_T(1-T(z)^*)^{-1}(1-T(z))^{-1}\Delta_T$.  Then the
curvature $K(T)$ of Arveson is defined by \cite{curvature}
\begin{equation}\label{E:curvaturefunction}
K(T)=\int_{S_d} \! \lim_{r\to 1}(1-r^2)tr(F(rz)) dz,
\end{equation}
where $dz$ is  the normalized Lebesgue measure on the unit sphere
$S_d$ in $\cc^d$.

\medskip

Results of Greene-Richter-Sundberg \cite{Greene} shows that the
curvature is known to be an integer. But it is not clear   how it
can be computed in terms of the actions of the operators $T_i$
(see Arveson's \cite{Arveson-psum} for more comments). Here we
shall show that the Samuel multiplicity provides such a way to
compute the curvature, and that it is obviously an integer.

\medskip

On the other hand, given a pure d-contraction with finite defect
rank, \emph{i.e.}, such that $I-T_1T_1^* - \cdots - T_dT_d^*$ is
finite rank, it follows that $H/(T_1H+\cdots +T_dH)$ is finite
dimensional. Hence, by considering $H$ as a Hilbert module over
$\cc[z_1, \cdots, z_d]$, the Samuel multiplicity with respect to
$I=(z_1, \cdots, z_d)$
\begin{equation}
e(I, H)=d! \ \lim_{k\to \infty} \frac{dim \ H/I^kH}{k^d}
\end{equation}
 exists, and is an integer.

The main result of this subsection is
\begin{theorem}\label{T:cur}
Let $T=(T_1, \cdots, T_d)$ be a pure d-contraction with finite
defect rank, acting on a Hilbert space $H$. Regard $H$ as a
Hilbert module over $\cc[z_1, \cdots, z_d]$, and define its Samuel
multiplicity $e(I, H)$ with respect to $I=(z_1, \cdots, z_d)$.

Then the curvature of $T$ is always equal to the Samuel
multiplicity, that is, $K(T)=e(I, H)$.
\end{theorem}
\textbf{Proof.} In the theory of pure d-contractions with finite
defect rank, a significant reduction, due to Arveson
\cite{subalgebra}, is that all these tuples can be realized as the
compressions of the tuple of multiplications by coordinate
functions onto coinvariant subspaces of the vector-valued
symmetric Fock space $H^2_d$ over the unit ball in $\cc^d$. Our
proof relies on this reduction.

Fix an invariant subspace $\submo \subset H^2_d\otimes \cc^N$ $(N
\in \nn)$ of the $\cc^N$-valued symmetric Fock space in $d$
variables. Let $T=(T_1, \cdots, T_d)$ acting on $H=\submoperp$ be
the compression  of the multiplication tuple $M_z=(M_{z_1},
\cdots, M_{z_d})$ acting on $H^2_d\otimes \cc^N$ onto
$\submoperp$. Then results in \cite{Greene} implies that
 $K(T)+f.d.(\submo)=N$. By the above discussions on
 $d! \lim_{k\to \infty}\frac{dim \ P_k {\submo}}{k^d}$, and Theorem
 \ref{T:identify}, we can complete the proof. $\quad \square$

 \bigskip

In \cite{Dirac} Arveson   pointed out that ``\emph{unlike the
Fredholm index, the curvature invariant is not known to be
invariant under similarity.}" Here as a quick application of
Theorem \ref{T:cur}, we show that it is.

\begin{corollary}\label{C:sim}
If $T=(T_1, \cdots, T_d)$ acting on $H$, and $S=(S_1, \cdots,
S_d)$ acting on $K$ are two d-contractions with finite defect
rank, and $T$ is similar to $S$; that is, there exists an
invertible operator $X\in B(H, K)$ such that $S_iX=XT_i$ for all
$i$, then $K(T)=K(S)$.
\end{corollary}
Proof. By arguments similar to the proof of Theorem 2, part (1) in
\cite{HP}, we know that the Hilbert polynomial $\phi_H(k)$, hence
the Samuel multiplicity $e(I, T)$, is invariant under similarity.
By Theorem \ref{T:main}, so is $K(T)$. $\quad \square$

\bigskip

Note that $\varphi_{\submo}(k+1)=dim \ P_k\submo=rank \
P_kP_{\submo}$, where $P_{\submo}$ is the projection onto
$\submo$. Also it is known \cite{HP}, \cite{Greene} that $d!
\lim_{k\to \infty} \frac{tr(P_kP_{\submo})}{k^d} $ is equal to the
fibre dimension $f.d.(\submo)$. Thus it follows that
\begin{corollary}
For any invariant subspace $\submo \subset H^2_d\otimes \cc^N$,
\begin{equation*}
d! \lim_{k\to \infty} \frac{tr(P_kP_{\submo})}{k^d}= d! \lim_{k\to
\infty} \frac{rank(P_kP_{\submo})}{k^d}.
\end{equation*}
\end{corollary}
Note that for two projections $P$ and $Q$, $tr(PQ)=rank(PQ)$ if
and only if $PQ$ is a projection, which means that the range of
$P$ naturally splits into a direct sum with respect to the range
and the kernel of $Q$. It suggests that  an invariant subspace
$\submo  \subset H^2_d\otimes \cc^N$ is asymptotically spliting
with respect to polynomials of large degrees.

We do not know how to give a direct proof of the above corollary.

\bigskip

Finally, in this subsection, we give an interesting interpretation
of the curvature invariant in terms of dilation theory. For any
invariant subspace $\submo \subset H^2_d\otimes \cc^N$, let
$H=\submoperp$, then  we have
\begin{equation*}
\phi_H(k)=dim \ H/I^kH = dim \ \submoperp\cap F_{k-1}.
\end{equation*}
So in terms of dilation theory, the
curvature $K(T)$ measures, in some sense, how many polynomials
$H=\submoperp$ contains.

\subsection{Proof of Theorem
\ref{T:identify}}\label{subS:identifyproof} The key to the proof
is to introduce an auxiliary invariant on $\submo$, denoted by
$\varepsilon_{\submo}$, and show that it is equal to the two
invariants appearing in Theorem \ref{T:identify}, respectively.
\begin{definition}\label{D:e}
For an invariant subspace $\submo \subset H^2_d\otimes \cc^N$,
define $\varepsilon_{\submo}(=\varepsilon)$ to be the maximal
dimension of a subspace $\mathcal{E}$ of $\cc^N$ with the
following property: there exists an orthonormal basis $e_1,
\cdots, e_{\varepsilon}$ of $\mathcal{E}$ and $h^1, \cdots,
h^{\varepsilon} \in \submo$ such that
\begin{equation*}
P_{H^2_d\otimes \mathcal{E}}h^i (\ne 0) \in H_2^d\otimes e_i,
\qquad i=1, \cdots, \varepsilon.
\end{equation*}
When $\mathcal{E}$ has the above property we say that   $\submo$
occupies $H^2_d\otimes \mathcal{E}$ in $H^2_d\otimes \cc^N$.
\end{definition}
\textbf{Remark}  It should be pointed out that $\varepsilon$ is
different from the smallest dimension of a subspace
$\mathcal{S}\subset \cc^N$, such that $\submo \subset H^2_d
\otimes \mathcal{S}$. A simple example can illustrate the
difference: take  $N=2$; for any two functions $f, g \in H^2_d$,
let $\submo$ be the invariant subspace generated by the
$\cc^2$-valued function $(f, g)$; then $\varepsilon=1$, while the
above $\mathcal{S}$ can be chosen to be one-dimensional if and
only if $f$ and $g$ differ by a scalar multiple.

\bigskip

 Next we give an elementary lemma which shows that the above definition is
independent of the basis $e_i$. It might be tempting to guess that
this conclusion has nothing to do with Nevanlinna-Pick reproducing
kernels, but just depends on the weighted shift structure. But
that is not the case. In fact, one can show that, the conclusion
is not true for the Bergman space over the unit disc, using the
fact that, two invariant subspaces of the scalar valued Bergman
space may have a positive angle.

\begin{lemma}\label{L:occupy}
If $\submo$ occupies $H^2_d\otimes \mathcal{E}$ for some
$\mathcal{E}\subset \cc^N$, then for any vector $e (\ne 0)\in
\mathcal{E}$, there exists an element $h\in \submo$ such that
$P_{H^2_d\otimes \mathcal{E}}h (\ne 0) \in H_2^d\otimes e$.
\end{lemma}
Proof. Fix an orthonormal basis $e_1, \cdots, e_N$ for $\cc^N$.
Then we write any element $f\in H^2_d\otimes \cc^N$ as $f=(f_1,
\cdots, f_N)$ with respect to the  decomposition $\oplus_{i=1}^{N}
H^2_d\otimes e_i$.  We say that $f$ has multiplier entries if each
$f_i$ is a multiplier on $H^2_d$. It is easy to see that the
condition that $f$ has multiplier entries is independent of the
choice of the basis $e_i$. Because $H^2_d$ admits a
Nevanlinna-Pick reproducing kernel, a theorem of McCullough and
Trent \cite{McCullough} implies that any invariant subspace
$\submo \subset H^2_d\otimes \cc^N$ is generated by elements with
multiplier entries. In particular, for any element $h (\ne 0) \in
H^2_d\otimes \cc^N$ the invariant subspace generated by $h=(h_1,
\cdots, h_N)$ contains a nonzero $f=(f_1, \cdots, f_N)$ which has
multiplier entries. Since
  $f=gh$ for some holomorphic $g$,   we know that $f_i=0$ if and only if
$h_i=0$. So in Definition \ref{D:e}, we can assume  that all $h^i$
have multiplier entries. In particular, let $P_{H^2_d\otimes
\mathcal{E}}h^i = g_i \otimes e_i$, where $g_i$ is a multiplier.
Then for any $e=c_1e_1+\cdots+c_{\varepsilon}e_{\varepsilon} \in
\mathcal{E}$ $(c_i\in \cc)$, let
\begin{equation*}
h=c_1g_2\cdots g_{\varepsilon}h^1+\cdots+c_{\varepsilon}g_1\cdots
g_{\varepsilon-1}h^{\varepsilon},
\end{equation*}
 then $P_{H^2_d\otimes \mathcal{E}}h =g_1\cdots g_{\varepsilon}\otimes e \in H_2^d\otimes e$. $\quad
 \Box$

\bigskip

For convenience, let $\delta=f.d.(\submo)$. The next lemma gives
the first, and easier, equality needed to prove Theorem
\ref{T:identify}.
\begin{lemma}\label{L:epsilon-fibre}
For any   $\submo \subset H^2_d\otimes \cc^N$, we have
$\delta=\varepsilon$.
\end{lemma}
Proof. It is obvious that $\delta\geq\varepsilon$. For the other
direction, choose $\delta$ elements $h^1, \cdots, h^{\delta}$ from
$\submo$, all with multiplier entries,  such that at some point
$\lambda \in B_d$ the $\delta$ vectors $h^1(\lambda), \cdots,
h^{\delta}(\lambda)$ are linearly independent in $\cc^N$. Now we
write $h^i=(h_1^i, \cdots, h_N^i)$ with respect to some basis
$e_1, \cdots, e_N$ of $\cc^N$, and assume that   the determinant
of the $\delta \times \delta$ matrix $\Theta=(h_j^i)_{i,
j=1}^{\delta}$ is nonzero. Note that the determinant $det(\Theta)$
is still a multiplier on $H^2_d$. Recall that the inverse   matrix
of $\Theta$ is given by $\frac{1}{det(\Theta)} (A_{i, j})_{i,
j=1}^{\delta}$, where $A_{i, j}$ is the $(\delta-1)\times
(\delta-1)$ minor of $\Theta$ associated with $h_j^i$. It follows
that $(h_j^i)(A_{i, j})=det(\Theta)\cdot \textbf{1}_{\delta}$ at
the level of matrix multiplication. For $j=1, \cdots, \delta$, if
we set $g^j=\sum_{i=1}^{\delta}A_{i, j}h^i$, then $P_{H^2_d\otimes
F}g^j=det(\Theta)e_j$, where $F=e_1\cc+\cdots+e_{\delta}\cc$. So
$\submo$ occupies $H^2_d\otimes F$. $\quad \Box$

\bigskip

The other equality needed to prove Theorem \ref{T:identify} is
given by Lemma \ref{L:last}.

\bigskip

For any element $f\in H^2_d \otimes \cc^N$, we expand $f$ into
homogeneous terms $f=f_c+f_{c+1}+\cdots$, where $f_i$ is a
homogeneous polynomial of degree $i$, and $f_c\ne 0$. Then we call
$c$ the order of $f$ at the origin, denoted by $ord(f)=c$. Let
$\mathcal{P}_k$ denote the set of all (scalar-valued) polynomials
with degrees at most $k$. (Note that we have used a different
notation $F_k$ to denote the vector-valued polynomials.) When
$N=1$, $\submo=H^2_d$,
we have $\varphi_{\submo}(k)=dim \ \mathcal{P}_{k-1} = \left(\begin{array}{c} d+k-1 \\
d\end{array}\right)$ and $\sigma=d! \lim_{k\to
\infty}\frac{\varphi_{\submo}(k)}{k^d}=1$.

\begin{lemma}\label{L:last}
For any  $\submo \subset H^2_d\otimes \cc^N$, we have $\sigma=d!
\lim_{k\to \infty}\frac{\varphi_{\submo}(k)}{k^d}=\varepsilon$.
\end{lemma}
Proof.  We first show that $\sigma=d! \lim_{k\to
\infty}\frac{\varphi_{\submo}(k)}{k^d}=1$ when  $\submo=[f]$, the
invariant subspace generated by a single element $f\in
H^2_d\otimes \cc^N$.   Note that
\begin{eqnarray*}
 dim \ P_{k+ord(f)}[f]  & \geq &
dim \ P_{k+ord(f)}(span\{\mathcal{P}_kf \})) \\
& = & dim \ \mathcal{P}_k \\
&  = & \left(\begin{array}{c} d+k \\
d\end{array}\right).
\end{eqnarray*}
 The first equality holds because  if $P_{k+c}(pf)=0$ for some
 polynomial with degree at most $k$, then we must have $p=0$.
  Now it is easy to see $\sigma=1$.

  Choose a  basis $e_1, \cdots, e_N$ for $\cc^N$, and write any
  $f=(f_1, \cdots, f_N)\in H^2_d\otimes \cc^N$ relative to the
  basis. Assume that $\submo$ occupies  $H^2_d\otimes \mathcal{E}$,
  where $\mathcal{E}=e_1\cc+\cdots+e_{\varepsilon}\cc$. Choose
  $h^i\in \submo$ such that $P_{H^2_d\otimes
  \mathcal{E}}h^i=h^i_i\otimes e_i ( \ne 0)$.
  Let $\submo'=[h^1_1\otimes e_1, \cdots,
  h^{\varepsilon}_{\varepsilon}\otimes e_{\varepsilon}]$, which naturally splits into
  the direct sum of $\varepsilon$ many singly generated invariant subspaces.
  Note that $\submo'\subset P_{H^2_d\otimes
  \mathcal{E}}(\submo)$. Now
    \begin{eqnarray*}
  dim \ P_k(\submo) & \geq & dim \ P_{H^2_d\otimes
  \mathcal{E}}P_k(\submo) \\
  & = & dim \ P_kP_{H^2_d\otimes
  \mathcal{E}}(\submo)\\
  & \geq & dim \ P_k(\submo') \\
  & = & \sum_{i=1}^{\varepsilon}dim \ P_k([h^i_i\otimes e_i]).
  \end{eqnarray*}
  It follows that $\sigma\geq \varepsilon$.

  For the other direction we first recall that for any $\lambda\in B_d$,
  $dim(\submo(\lambda))\leq \delta$, which is equal to $\varepsilon$ by Lemma
  \ref{L:epsilon-fibre}. Hence for any $f=(f_1, \cdots, f_N)\in
  \submo$, the determinant of the $(\varepsilon+1)\times
  (\varepsilon+1)$ matrix
  \begin{equation*}
  \left(\begin{array}{ccccc}h^1_1 & 0 & \cdots & 0 & f_1 \\
  0 & h^2_2 & \cdots & 0 & f_2 \\
\cdots &\cdots &\cdots &\cdots &\cdots \\
0 & 0  & \cdots & h^{\varepsilon}_{\varepsilon} & f_{\varepsilon}
\\
h^1_i & h^2_i & \cdots & h^{\varepsilon}_i & f_i \end{array}
\right)
\end{equation*}
is identically zero for any fixed $i=\varepsilon+1, \cdots, N$. It
follows that
\begin{equation}\label{E:key}
g_1f_1+g_2f_2+\cdots+g_{\varepsilon}f_{\varepsilon}+g_if_i=0,
\end{equation}
where $g_j=h^1_1\cdots h^{j-1}_{j-1}\cdot h^j_i \cdot
h^{j+1}_{j+1}\cdots h^{\varepsilon}_{\varepsilon}$ for $j=1,
\cdots, \varepsilon$, and $g_i=h^1_1\cdots
h^{\varepsilon}_{\varepsilon}$. In particular, $g_i$ is nonzero,
and independent of $i$.

Now for   $k\in \nn$, we consider the natural map   $J_k:
P_k(\submo)\to P_{H^2_d\otimes\mathcal{E}}P_k(\submo)$. If
$\xi=P_k(f) \in ker(J_k)$, that is,
\begin{equation*}
P_k(f_1\otimes e_1)=\cdots=P_k(f_{\varepsilon}\otimes
e_{\varepsilon})=0,
\end{equation*}
 then by Eq. \ref{E:key}
 \begin{equation*}
P_k(g_1f_1+\cdots+g_{\varepsilon}f_{\varepsilon})=-P_k(g_if_i)=0.
\end{equation*}
 Note that $ord(g_i)$  is independent of $i$. Now by looking at the lowest
degree term in $g_if_i$, we conclude that $P_{k-ord(g_i)}f_i=0$
and  $i=\varepsilon+1, \cdots, N$. This implies that the kernel
$ker(J_k)$ is contained in the range of
$P_{H^2_d\otimes\mathcal{E}^{\perp}}(P_k-P_{k-ord(g_i)})$ whose
rank is
\begin{displaymath}
(N-\varepsilon)(\left(\begin{array}{c} d+k \\
d\end{array}\right)-\left(\begin{array}{c} d+k-ord(g_i) \\
d\end{array}\right)),
\end{displaymath}
 which is a polynomial of degree $d-1$. Hence
\begin{equation*}
\sigma=d! \lim_{k\to \infty}\frac{dim \ P_k(\submo)}{k^d} = d!
\lim_{k\to \infty}\frac{dim \
P_{H^2_d\otimes\mathcal{E}}P_k(\submo) }{k^d}\leq \varepsilon.
\qquad \Box
\end{equation*}
\section{Concluding remarks}
 Since part of our purpose is to develop a
multivariable Fredholm theory,  we give a different proof of the
  index formula $index(M_z, M_w)=e_0$,  appearing in Theorem \ref{T:subFock}.
   But we are not able to
find a different proof  for the characterization of the
Fredholmness of $(M_z, M_w)$.

We have to assume that $(M_z,M_w)$ is Fredholm.
Then by step 3, and 6 in the proof of Theorem \ref{T:subFock}, we
have
\begin{equation}
index(T_w)=-\lim_{k \to \infty} \frac{dim \ \submo/(z\submo +w^k
\submo)}{k}.
\end{equation}

By Proposition \ref{P:reduction-mul}
\begin{equation*}
\lim_{k \to \infty} \frac{dim \ \submo/(z\submo +w^k \submo)}{k}
\ge e(I, \submo) =e_0,
\end{equation*}
where $I=(z, w)\subset \cc[z, w]$ is the maximal ideal at the
origin.

By the result of Gleason-Richter-Sundberg \cite{Gleason}, we know
that, if $\submo$ is Fredholm, then
\begin{equation*}
index(M_z, M_w) = f.d. (\submo).
\end{equation*}

By Theorem \ref{T:identify}, $f.d.(\submo)=d! \ \lim_{k \to
\infty} \frac{dim \ P_k \submo}{k^d}$.

Note that the kernel of the following natural surjective map
\begin{displaymath}
\submo \to P_k\submo
\end{displaymath}
contains $I^{k+1}\submo$, and hence it factors through
\begin{displaymath}
\submo / I^{k+1} \submo \to P_k\submo.
\end{displaymath}
It follows that
\begin{displaymath}
e(I, \submo)= d! \ \lim_{k \to \infty} \frac{dim \ \submo /
I^{k+1} \submo}{k^d} \ge d! \ \lim_{k \to \infty} \frac{dim \ P_k
\submo}{k^d} = f.d.(\submo).
\end{displaymath}
Now we have another proof of the index formula in Theorem
\ref{T:subFock}.

\bigskip

 Before we conclude the paper, we make some comments on further
 studies. Hopefully, this will spur the interests of some readers.

 So far, we have several
numerical invariants defined on $\submo \subset H^2 \otimes
\cc^N$, appearing in \cite{HP}, \cite{Fred}, \cite{Dirichlet},
\cite{Hardy},  and this paper:

\medskip

 1. $index(M_z, M_w)$;

 2. $e(I, \submo)$;

 3. $\lim_{k \to \infty}
 \frac{dim \ \submo/ (z \submo + w^k \submo)}{  k}$;

 4. $dim \ \submo/[(z-\lambda)\submo +
 (w-\mu)\submo]$ for almost all $(\lambda, \mu)$ in the ball;

 5. $ d! \ \lim_{k \to \infty} \frac{rank \ P_k P_{\submo}}{k^d} $,
here $P_{\submo}$ is the orthogonal projection onto $\submo$;

6. $ d! \ \lim_{k \to \infty} \frac{trace \ P_k P_{\submo}}{k^d}$;

7. $f.d.(\submo)$;

8. the $\varepsilon$ invariant defined in
\ref{subS:identifyproof}.

\bigskip

Among them the last five are always well-defined, and equal;
moreover, when $dim \ \submo/I \submo < \infty$, the first three
are defined, and are equal to the last five.

 It seems that the above equalities, or the failure of these
 equalities, together with the additivity of Hilbert polynomials as in
 Theorem \ref{T:add} and \cite{Dirichlet},
  can serve as test problems when trying to develop
 operator theory in several variables or over different spaces.
  It
 is interesting to observe that, when one looks at the Bergman
 space over various domains, almost all the above equalities can fail to be true.
  This might pose a host of problems. For instance, is
 the invariant in the above (6) always well defined for \emph{any} invariant
 subspace of the Bergman space $L^2_a(\dd)$ over the unit disc? A
 positive answer would
 provide a numerical invariant, lying between $0$ and $1$,
  measuring the size of the invariant subspaces of
 $L^2_a(\dd)$;
in fact, if we assume its existence, then examples show that it
takes all the values in the interval $(0, 1]$.



\noindent
Department of Mathematics\\
University of Alabama\\
Tuscaloosa, AL 35487 \\
xfang@bama.ua.edu

\bigno \emph{current address:}\\
Department of Mathematics\\
Kansas State University\\
Manhattan, KS 66502\\
xfang@math.ksu.edu
\end{document}